\documentclass[
final
]{dmtcs-episciences}
\pdfoutput=1

\usepackage[utf8]{inputenc}
\usepackage{amssymb}
\usepackage{amsmath}
\usepackage{graphicx}

\usepackage{tikz}
\usetikzlibrary{calc}
\usetikzlibrary{decorations.pathreplacing}

\usetikzlibrary{external}
\newif\ifcomment
\commentfalse
\commenttrue

\newif\ifcochain
\cochainfalse
\cochaintrue

\newif\iflong
\newif\ifshort

\longtrue

\iflong
\else
\shorttrue
\fi
\usepackage{amsthm}

\usepackage{mathtools}
\usepackage{xcolor}
\usepackage{hyperref}
\usepackage{cleveref}
\usepackage{ECmacros}
\usetikzlibrary{calc}
\usepackage{ifthen}
\usepackage{nicefrac}
\usepackage{enumerate}
\usepackage[absolute]{textpos}
\usepackage{caption}
\usepackage[labelformat=simple]{subcaption}
\usepackage[T1]{fontenc} % so that Ł looks better in `Łukasz'

\newtheorem{theorem}{Theorem}
\newtheorem{lemma}[theorem]{Lemma}
\newtheorem{corollary}[theorem]{Corollary}

\newtheorem{conjecture}{Conjecture}

\author[M. Bonamy, Ł. Bożyk, A. Grzesik, M. Hatzel, T. Masařík, J. Novotná, and K. Okrasa]{Marthe Bonamy
\affiliationmark{1}
\and \L{}ukasz Bo\.{z}yk
\affiliationmark{2}
\and Andrzej Grzesik
\affiliationmark{3}
\and Meike Hatzel
\affiliationmark{4}\thanks{M.~Hatzel was supported by a fellowship within the IFI programme of the German Academic Exchange Service (DAAD).}\\
\and Tom\'a\v{s} Masa\v{r}\'ik
\affiliationmark{2,5}\thanks{T.~Masařík was supported by a postdoctoral fellowship at SFU (NSERC grants R611450 and R611368).}
\and Jana Novotn\'{a}
\affiliationmark{2}
\and Karolina Okrasa
\affiliationmark{2,6}
}
\title{Tuza's Conjecture for Threshold Graphs
\thanks{
This research has received funding from the European Research Council (ERC) under the European Union's Horizon 2020 research and innovation programme, grant agreements No.~714704 (Ł.~Bożyk, T.~Masařík, J.~Novotná, and K.~Okrasa), No.~648509 (A.~Grzesik), and No.~648527 (M.~Hatzel).
}
\thanks{An extended abstract of this manuscript has been accepted at EUROCOMB 2021~\cite{TuzaEurocomb}.}
}
\affiliation{
CNRS, Universit\'e de Bordeaux, France\\
Institute of Informatics,
University of Warsaw, Poland\\
Jagiellonian University, Kraków, Poland\\
National Institute of Informatics, Tokyo, Japan\\
Simon Fraser University, Burnaby, BC, Canada\\
Warsaw University of Technology, Poland}


\received{2021-07-09}

\revised{2022-03-16}

\accepted{2022-06-03}

\begin{document}

\publicationdetails{24}{2022}{1}{24}{7660}
\maketitle

\begin{textblock}{20}(0, 12.5)
\includegraphics[width=40px]{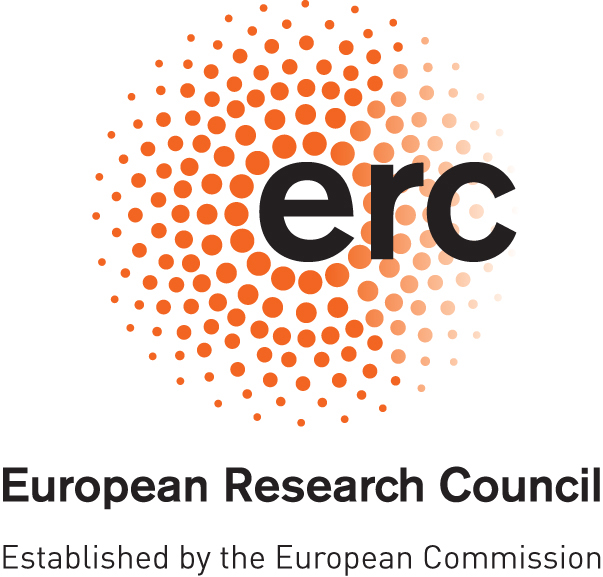}%
\end{textblock}
\begin{textblock}{20}(-0.25, 12.9)
\includegraphics[width=60px]{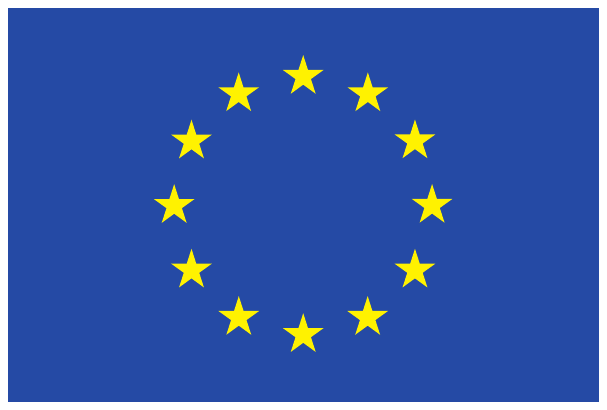}%
\end{textblock}

\begin{abstract}
Tuza famously conjectured in 1981 that in a graph without $k+1$ edge-disjoint triangles, it suffices to delete at most $2k$ edges to obtain a triangle-free graph.
The conjecture holds for graphs with small treewidth or small maximum average degree, including planar graphs. However, for dense graphs that are neither cliques nor 4-colourable, only asymptotic results are known.
Here, we confirm the conjecture for threshold graphs, i.e.\@ graphs that are both split graphs and cographs, and for co-chain graphs with both sides of the same size divisible by $4$.
\end{abstract}

\renewcommand{\TH}[1]{\ensuremath{\tau(#1)}}
\newcommand{\TP}[1]{\ensuremath{\mu(#1)}}
\newcommand{\Kt}{K_{\operatorname{top}}}
\newcommand{\Kb}{K_{\operatorname{bot}}}
\newcommand{\Kto}{K^{\operatorname{top}}_1{}}
\newcommand{\Kbo}{K^{\operatorname{bot}}_1{}}
\newcommand{\Ktt}{K^{\operatorname{top}}_2{}}
\newcommand{\Kbt}{K^{\operatorname{bot}}_2{}}

\newcommand{\THset}{\mathsf{TH}}
\newcommand{\TPset}{\mathsf{TP}}

\section{Introduction}

If we can ``pack'' at most $k$ disjoint objects of some type  in a given graph, how many elements do we need to ``cover'' all appearances of such an object in the graph? Erd\H{o}s and P\'osa famously proved that if a graph contains at most $k$ pairwise vertex-disjoint cycles, then there is a set of at most $f(k)$ vertices that intersects every cycle~\cite{erdosposa}.
While the exact best value of function~$f$ is yet unknown, the asymptotic behaviour was recently determined to be $f(k)=\Theta(k \log k)$~\cite{van2019tight}. 

In this paper, we focus on edge-disjoint triangles; we refer the interested reader to~\cite{JF} for a dynamic survey on other objects. 
For a graph $G$, we call every family of pairwise edge-disjoint triangles a \emph{triangle packing}, and every subset of edges intersecting all triangles in $G$ a \emph{triangle hitting}. We denote by $\TP{G}$ the maximum size of a triangle packing in $G$, and by $\TH{G}$ the minimum size of a triangle hitting in $G$.
Trivially, there is a set of at most $3\TP{G}$ edges that intersect every triangle. We are concerned with improving that bound, following Tuza's conjecture from 1981.

\begin{conjecture}[Tuza \cite{tuza1981conj}]\label{conj:Tuza}
    For any graph $G$ it holds $\TH{G} \leq 2 \TP{G}$.
\end{conjecture}

\Cref{conj:Tuza}, if true, is tight for $K_4$ and $K_5$.
Gluing together copies of $K_4$ and $K_5$ along vertices, it is easy to build an infinite family of connected graphs for which \cref{conj:Tuza} is tight.
However, for larger cliques, it is known that the ratio $\TH{K_p}/\TP{K_p}$ tends to $3/2$ as $p$ increases \cite{feder2012packing}.
In addition, Haxell and R\"{o}dl~\cite{HaxellRodl} proved that $\TH{G} \leq 2 \TP{G}+o(|V(G)|^2)$ for any graph $G$, meaning \cref{conj:Tuza} is asymptotically true when $\TH{G}$ is quadratic with respect to $|V(G)|$.
Those seem to indicate that \cref{conj:Tuza} should be easier for dense graphs than for sparse graphs. 
Conversely, it is asymptotically tight in some classes of dense graphs~\cite{tightfordensegraphs}.
If we focus on \emph{hereditary graph} classes (i.e.\@ classes that contain every induced subgraph of a graph in the class), the conjecture has only been confirmed for a few graph classes. Those classes include most notably graphs of treewidth at most~$6$~\cite{botler2021tuza}, 4-colourable graphs~\cite{4colorable}, and graphs with maximum average degree less than $7$~\cite{puleo2015tuza}.

A good candidate for an interesting dense hereditary graph class is the class of \emph{split graphs}, i.e.\@ graphs whose vertex set can be partitioned into two sets: one that induces a clique, the other inducing an independent set.
However, \cref{conj:Tuza} remains a real challenge even when restricted to split graphs. 
Another good candidate for an interesting dense hereditary graph class is the class of \emph{cographs}, i.e.\@ graphs with no induced path on four vertices. As an initial step, we focus on graphs that are both split graphs and cographs, i.e.\@ \emph{threshold} graphs. 
While this may seem like a small step, it is arguably the first dense hereditary superclass of cliques where the conjecture is confirmed. 

\begin{theorem}\label{th:main}
    If $G$ is a threshold graph, then $\TH{G} \leq 2 \TP{G}$.
\end{theorem}

In the latter part of the paper, we show that similar tools with more involved analysis can be used to verify \cref{conj:Tuza} also for specific co-chain graphs.
A graph $G$ is a \emph{co-chain graph} (or sometimes alternatively called \emph{co-difference graph}) if its vertex set can be partitioned into two sets $K_1$ and $K_2$ such that $G[K_1]$ and $G[K_2]$ are cliques and there is an ordering $c_1,\dots,c_n$ on the vertices of $K_1$ and an ordering $d_1,\dots,d_m$ on the vertices of $K_2$ with $N[c_{i+1}] \subseteq N[c_i]$ for all $1 \leq i < n$ and $N[d_i] \subseteq N[d_{i+1}]$ for all $1 \leq i < m$.
We call $\Brace{K_1,K_2}$ a \emph{co-chain representation} of $G$.
We say that $G$ is an \emph{even balanced} co-chain graph if additionally $K_1$ and $K_2$ are of the same size that is divisible by four.

\begin{theorem}\label{thm:cochain}
    If $G$ is an even balanced co-chain graph, then $\TH{G} \leq 2 \TP{G}$.
\end{theorem}
 
\Cref{thm:cochain} can be seen as a very first step towards attacking \cref{conj:Tuza} on (mixed) unit interval graphs as those graphs can be modelled as a \emph{concatenation} of co-chain graphs.
That is, vertices of graph $G$ are partitioned into $r$ cliques $C_1,\ldots,C_r$ where each $\Brace{C_i,C_{i+1}}$ induce a co-chain graph and $G$ contains no other edges; see~\cite{unit,mixed} for more details.
The simplest object for further study might be a $k$-path, which can be viewed as a concatenation of well-structured same-sized co-chain graphs.

Finally, it is worth mentioning that \Cref{conj:Tuza} is known to hold as soon as we consider \emph{multi}-packing~\cite{chalermsook2020multi}, and in particular it holds in its fractional relaxation~\cite{Krivelevich1995}. Another angle of attack consists of lowering the bound of $3$ step by step for all graphs.
The best, and in fact only, such bound is slightly under $2.87$~\cite{haxell1999packing}.

\subsection{Preliminaries}

All graphs in this paper are undirected and simple.
Let $G = (V,E)$ be a graph.
By the \emph{size} of a graph $G$ (alt.~$\Abs{G}$), we always mean the number of its vertices.
For all $v \in V$ the set $N(v) \coloneqq \Set{u \mid \Set{u,v}\in E}$ is called the \emph{neighbourhood} of $v$ and $N[v] \coloneqq N(v) \cup \Set{v}$ is its \emph{closed neighbourhood}.
A \emph{matching} in $G$ is a set of edges $M \subseteq E$ such that every vertex of $G$ is incident to at most one edge of $M$.
A vertex $v\in V$ is \emph{complete} to $A\subseteq V, v\notin A$ if $v$ is adjacent to all vertices in $A$. 
Disjoint sets $A,B\subseteq V$ are \emph{complete} to each other if $E$ contains all edges between $A$ and $B$. 
Any omitted definitions can be found in the book by Diestel \cite{diestelbook}. 

\begin{figure}[htb]
        \centering
        \begin{tikzpicture}[scale=\ifshort0.6\fi\iflong0.72\fi]
                \node at (-4,0) {\begin{tikzpicture}
		\tikzstyle{vertex}=[shape=circle, fill=black, draw, inner sep=.6mm]
		\definecolor{c-1}{rgb}{0.33, 0.42, 0.18}
		\definecolor{c-3}{rgb}{0.55, 0.71, 0.0}
		\definecolor{c-4}{rgb}{0.81, 0.09, 0.13}
		\definecolor{c-2}{rgb}{0.0, 0.5, 1.0}
		\definecolor{c-0}{rgb}{0.41, 0.16, 0.38}
		\definecolor{c-5}{rgb}{1.0, 0.51, 0.26}
		\tikzset{matching/.style={line width=2pt,color=#1}}
		\node (center) at (0,0) {};
		\node[vertex] (v-0) at (0,0) {};
		
		\foreach \x in {1,...,5}
		{
			\node[vertex] (v-\x) at ($(center)+({360/5*\x}:1.4)$) {};
			\node (v-\x-l) at ($(v-\x)+({360/5*\x}:0.35)$) {$\x$};
		}
		\foreach \x in {1,...,5}
		{	
			%M_x
			
			\draw[matching={c-\x}] (v-0) to (v-\x);
			\foreach \a in {1,...,5}
			{
				\foreach \b in {1,...,5}
				{
					\pgfmathtruncatemacro{\apb}{\a+\b};
					\pgfmathtruncatemacro{\y}{mod(\apb,5)};
					\pgfmathtruncatemacro{\xx}{mod(2*\x,5)};
					\ifthenelse{\NOT\a=\b \AND \(\y = \xx\)}{\draw[matching={c-\x}] (v-\a) to (v-\b);}{}
				}
			}
		}
		\node (v-0-l) at ($(v-0)+({325}:0.3)$) {$0$};
	\end{tikzpicture}};
	            \node at (4,0) {\begin{tikzpicture}
		\tikzstyle{vertex}=[shape=circle, fill=black, draw, inner sep=.6mm]
		\definecolor{c-1}{rgb}{0.33, 0.42, 0.18}
		\definecolor{c-3}{rgb}{0.55, 0.71, 0.0}
		\definecolor{c-4}{rgb}{0.81, 0.09, 0.13}
		\definecolor{c-2}{rgb}{0.0, 0.5, 1.0}
		\definecolor{c-0}{rgb}{0.41, 0.16, 0.38}
		\definecolor{c-5}{rgb}{1.0, 0.51, 0.26}
		\tikzset{matching/.style={line width=2pt,color=#1}}
		\node (center) at (0,0) {};
		
		\foreach \x in {1,...,5}
		{
			\node[vertex] (v-\x) at ($(center)+({360/5*\x}:1.4)$) {};
			\node (v-\x-l) at ($(v-\x)+({360/5*\x}:0.35)$) {$\x$};
		}
		\foreach \x in {1,...,5}
		{	
			%M_x
			
			\foreach \a in {1,...,5}
			{
				\foreach \b in {1,...,5}
				{
					\pgfmathtruncatemacro{\apb}{\a+\b};
					\pgfmathtruncatemacro{\y}{mod(\apb,5)};
					\pgfmathtruncatemacro{\xx}{mod(2*\x,5)};
					\ifthenelse{\NOT\a=\b \AND \(\y = \xx\)}{\draw[matching={c-\x}] (v-\a) to (v-\b);}{}
				}
			}
		}
	\end{tikzpicture}};
        \end{tikzpicture}
	    \caption{The decomposition of edges of a $6$-vertex clique 
	    into $5$ matchings and the corresponding decomposition of a $5$-vertex clique.}
        \label{fig:match-dec}
    \end{figure}
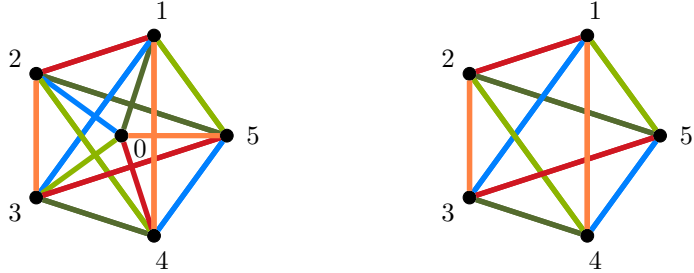

Let us first recall the following well-known property (chromatic index of a clique).
\begin{lemma}
    \label{thm:mat_decomp}
    The edge set of a clique $K$ 
    on $k$ vertices can be decomposed into $k$ edge disjoint maximal matchings for $k$ odd and $k-1$ edge disjoint maximal matchings for $k$ even.
\end{lemma}
\begin{proof}
    If $k$ is even, we may identify the vertices of $K$ with the set $\{0,1,\ldots,k-1\}$ and consider matchings
\begin{align*}
    M_i=\Set{\Set{0,i}}\cup\Set{\Set{a,b}\mid a\neq b,\ ab \neq 0,\ a+b\equiv 2i\pmod{k-1}}
\end{align*}
    for $1\leq i\leq k-1$.
    These matchings are edge disjoint and cover the entire edge set of $K$ (cf.\@ \cref{fig:match-dec}).
    Removing any vertex (along with all incident edges) yields a desired matching decomposition into $k-1$ matchings of the edge set of the clique of $k-1$ vertices.
\end{proof}
    
A graph $G=(V,E)$ is a \emph{star} if $V = \Set{c,s_1,\dots,s_k}$ and $E =  \Set{\Set{c,s_i}|1\leq i \leq k}$; the vertex $c$ is called the \emph{center vertex} of the star.
A graph $G$ is a \emph{complete split graph} if its vertex set can be partitioned into sets $K$ and $S$, such that $S$ is independent, $K$ induces a clique, and $K$ and $S$ are complete to each other.

The following lemma describes how to pack triangles in complete split graphs.
As it is very central to our proofs later, we include a proof here.

\begin{figure}[htb]
     \centering
     \begin{subfigure}[b]{0.49\textwidth}
         \centering
         \begin{tikzpicture}[xscale=\ifshort.5\fi\iflong.6\fi,yscale=\ifshort.4\fi\iflong.48\fi]

\def\rndcrn{5pt}
\def\crnshift{9pt}

\tikzstyle{vertex}=[circle,minimum size=4pt,inner sep=0pt,draw,fill]

%%% clique part coordinates

\coordinate (Ksw) at (0,-8);
\coordinate (Kse) at (2,-8);
\coordinate (Kne) at (2,0);

\coordinate (Klab) at (-1,-4);

\def\Kvx{1} % x position of vertices inside K

\coordinate (c1) at (\Kvx,-.5);
\coordinate (c2) at (\Kvx,-1.5);
\coordinate (c3) at (\Kvx,-2.5);
\coordinate (c4) at (\Kvx,-3.5);
\coordinate (ck1) at (\Kvx,-6.5);
\coordinate (ck) at (\Kvx,-7.5);

%%% stable part coordinates

\coordinate (Ssw) at (3,-7.5);
\coordinate (Snw) at (3,-.5);
\coordinate (Sne) at (5,-.5);

\coordinate (Slab) at (6,-4);

\def\Svx{4} % x position of vertices inside S

\coordinate (u1) at (\Svx,-1);
\coordinate (u2) at (\Svx,-2);
\coordinate (u3) at (\Svx,-3);
\coordinate (us1) at (\Svx,-6);
\coordinate (us) at (\Svx,-7);

%%% fill K

\draw[thick,fill=black!20,rounded corners=\rndcrn] (Ksw) rectangle (Kne);

%%% draw full join between S and K

\begin{scope}
\draw[fill=black!10] ([yshift=\crnshift]Ssw)--([yshift=\crnshift]Kse)--([yshift=-\crnshift]Kne)--([yshift=-\crnshift]Snw)--cycle;
\end{scope}

%%% draw parts

\draw[thick,rounded corners=\rndcrn] (Ksw) rectangle (Kne);
\draw[thick,rounded corners=\rndcrn] (Ssw) rectangle (Sne);

%%% draw triangles

\begin{scope}[opacity=.75]
\draw[color=blue,fill=blue!20] (u2) -- (c1) -- (c2)-- (u2) -- (c3) -- (c4) -- (u2) -- (ck1) -- (ck) -- cycle;
\draw[very thick,color=blue] (c1) edge (c2);
\draw[very thick,color=blue] (c3) edge (c4);
\draw[very thick,color=blue] (ck1) edge (ck);
\end{scope}

%%% label parts

\node at (Klab) {$K$};
\node at (Slab) {$S$};

%%% draw and label vertices

\node[vertex] at (c1) {};
\node[vertex] at (c2) {};
\node[vertex] at (c3) {};
\node[vertex] at (c4) {};
\node[vertex] at (ck1) {};
\node[vertex] at (ck) {};

\node at  ($ (c4)!2/5!(ck1) $) {$\vdots$};

\node[vertex] at (u1) {};
\node[vertex,color=blue] at (u2) {};
\node[vertex] at (u3) {};
\node[vertex] at (us1) {};
\node[vertex] at (us) {};

\node at ($ (u3)!2/5!(us1) $) {$\vdots$};

\end{tikzpicture}
         \label{fig:match-lift-left}
     \end{subfigure}
     \hfill
     \begin{subfigure}[b]{0.49\textwidth}
         \centering
         \begin{tikzpicture}[xscale=\ifshort.5\fi\iflong.6\fi,yscale=\ifshort.4\fi\iflong.48\fi]

\def\rndcrn{5pt}
\def\crnshift{9pt}

\tikzstyle{vertex}=[circle,minimum size=4pt,inner sep=0pt,draw,fill]

%%% clique part coordinates

\coordinate (Ksw) at (0,-8);
\coordinate (Kse) at (2,-8);
\coordinate (Kne) at (2,0);

\coordinate (Klab) at (-1,-4);

\def\Kvx{1} % x position of vertices inside K

\coordinate (c1) at (\Kvx,-.5);
\coordinate (c2) at (\Kvx,-1.5);
\coordinate (c3) at (\Kvx,-2.5);
\coordinate (c4) at (\Kvx,-3.5);
\coordinate (ck1) at (\Kvx,-6.5);
\coordinate (ck) at (\Kvx,-7.5);

%%% stable part coordinates

\coordinate (Ssw) at (3,-7.5);
\coordinate (Snw) at (3,-.5);
\coordinate (Sne) at (5,-.5);

\coordinate (Slab) at (6,-4);

\def\Svx{4} % x position of vertices inside S

\coordinate (u1) at (\Svx,-1);
\coordinate (u2) at (\Svx,-2);
\coordinate (u3) at (\Svx,-3);
\coordinate (us1) at (\Svx,-6);
\coordinate (us) at (\Svx,-7);

%%% fill K

\draw[thick,fill=black!20,rounded corners=\rndcrn] (Ksw) rectangle (Kne);

%%% draw full join between S and K

\begin{scope}
\draw[fill=black!10] ([yshift=\crnshift]Ssw)--([yshift=\crnshift]Kse)--([yshift=-\crnshift]Kne)--([yshift=-\crnshift]Snw)--cycle;
\end{scope}

%%% draw parts

\draw[thick,rounded corners=\rndcrn] (Ksw) rectangle (Kne);
\draw[thick,rounded corners=\rndcrn] (Ssw) rectangle (Sne);

%%% draw triangles

\begin{scope}[opacity=.75]
\draw[color=orange,fill=orange!20] (u3) -- (c4) to [out=-70,in=70] (ck1)-- (u3) -- (c2) -- (c3) -- cycle;
\end{scope}

\begin{scope}[opacity=.5]
\draw[color=orange,fill=orange!20] (u3) -- (c1) to [out=-70,in=70] (ck) -- cycle;
\end{scope}

\draw[color=orange,opacity=.75] (ck) --(u3) -- (c1);

\draw[very thick,color=orange] (c3) edge (c2);
\draw[very thick,color=orange] (c4) to [out=-70,in=70] (ck1);
\draw[very thick,color=orange] (c1) to [out=-70,in=70] (ck);

%%% label parts

\node at (Klab) {$K$};
\node at (Slab) {$S$};

%%% draw and label vertices

\node[vertex] at (c1) {};
\node[vertex] at (c2) {};
\node[vertex] at (c3) {};
\node[vertex] at (c4) {};
\node[vertex] at (ck1) {};
\node[vertex] at (ck) {};

\node at  ($ (c4)!2/5!(ck1) $) {$\vdots$};

\node[vertex] at (u1) {};
\node[vertex] at (u2) {};
\node[vertex,color=orange] at (u3) {};
\node[vertex] at (us1) {};
\node[vertex] at (us) {};

\node at ($ (u3)!2/5!(us1) $) {$\vdots$};

\end{tikzpicture}
         \label{fig:match-lift-right}
     \end{subfigure}
        \caption{Full joins of matchings in $K$ with vertices in $S$ as families of triangles.}
        \label{fig:match-lift}
\end{figure}
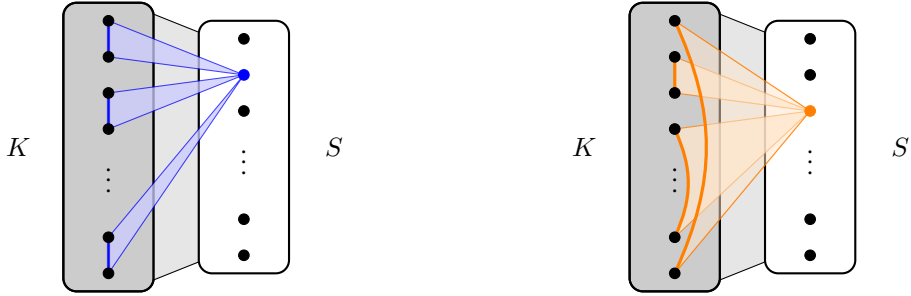

\begin{lemma}[\cite{feder2012packing}]
    \label{lem:clique_IS}
	Let $K$ be a clique, $S$ an independent set such that they are complete to each other and $\Abs{K}=\Abs{S}=k$. 
	Then we can find an (optimal) triangle packing $\TPset$ of size $\binom{k}{2}$ such that:
	\begin{enumerate} 
		\item It uses all edges from $K$ and each triangle in $\TPset$ contains exactly one edge from $K$.
		\item If $k$ is odd, the remaining edges (not used in $\TPset$) create a matching between $K$ and $S$, otherwise they create a star with its center vertex in $S$. Moreover, we can choose the unused matching and the center vertex of the unused star arbitrarily.
	\end{enumerate}
\end{lemma}
\begin{proof} 
    Consider a graph $G$ composed of a clique $K'$ complete to an independent set $S'$ with $\Abs{K'}=k$ and $\Abs{S'}=k-1$, where $k$ is even. 
    By \cref{thm:mat_decomp}, $K$ can be decomposed into $k-1$ edge disjoint (perfect) matchings of size $\nicefrac{k}{2}$.
    Each such matching fully joined to a different vertex in $S'$ yields a family of $\nicefrac{k}{2}$ edge disjoint triangles (see \cref{fig:match-lift}).
    The collection of all $k-1$ such joins is a decomposition of the entire edge set of $G$ into triangles.

    Removing any vertex $u$ from $K'$ yields a balanced graph with both sides of odd size, in which edges not packed into triangles (participating in triangles whose vertex $u$ got removed) create a matching between $K'-u$ and $S'$.
    On the other hand, by adding a single vertex $v$ to $S'$, we get a balanced graph with both sides of even size, in which unpacked edges form a star (with $v$ being its center vertex).
\end{proof}

\begin{corollary}
    \label{cor:clique_IS}
    Let $K$ be a clique and $S$ an independent set such that they are complete to each other. 
    \begin{enumerate}[(a)]
    \item If $\Abs{S}< \Abs{K}$, then we can find a triangle packing of size $\Abs{S}\cdot\Floor{\nicefrac{\Abs{K}}{2}}$. 
     \label{cor:clique_IS:between}
    \item If $\Abs{S} \geq \Abs{K}$, then we can find a triangle packing of size $\binom{\Abs{K}}{2}$.
    \label{cor:clique_IS:inside}
    \end{enumerate}
\end{corollary}
\begin{proof}
If $\Abs{S}< \Abs{K}$, we take arbitrary $\Abs{S}$ edge-disjoint maximal matchings in $K$ whose existence follows from \cref{thm:mat_decomp} and assign them to different vertices in $S$. 
The full join of each such pair consists of $\Floor{\nicefrac{\Abs{K}}{2}}$ edge-disjoint triangles.

If $\Abs{S} \geq\Abs{K}$, we can derive the statement from \cref{lem:clique_IS}: it is enough to take any $\Abs{K}$-element subset $S'$ of $S$.
\end{proof}

We say that we \emph{pack edges of $K$ with vertices of $S$} when we use triangle packings from \cref{cor:clique_IS}.
The following lemma describes tightly how many edge-disjoint triangles can be packed in a clique.
\begin{lemma}[\cite{feder2012packing}]
	\label{lem:clique}
	The optimal triangle packing for $K_n$ with $n=6x+i,0\le i\le5$ is $\nicefrac{\Brace{\binom{n}{2}-k}}{3}$ where $k$ is the number of not covered edges and
	\begin{itemize}
		\item $k=0$ for $i=1,3$,
		\item $k=4$ for $i=5$,
		\item $k=\frac{n}{2}$ for $i=0,2$,
		\item $k=\frac{n}{2}+1$ for $i=4$.
	\end{itemize}
\end{lemma}
Observe, that we can always hit all the triangles in a clique by leaving a bipartite graph with partitions of as equal size as possible and removing the rest.
Therefore, the optimal triangle hitting in a clique consists of at most half the edges.

\section{Threshold graphs}

A graph $G=(V,E)$ is a \emph{threshold graph} if its vertex set can be partitioned into two sets $K=\{c_1,\ldots,c_{k}\}$ and 
$S=\{u_1,\ldots,u_{s}\}$ such that $G[K]$ is a clique and $G[S]$ is an independent set in $G$, and
$N[c_{i+1}] \subseteq N[c_i]$ for all $1 \leq i < k$ and $N(u_i) \subseteq N(u_{i+1})$ for all $1 \leq i < s$.
We identify $K$ with the clique $G[K]$ and say $G=(K\cup S,E)$ is a threshold graph with given \emph{threshold representation} $(K,S)$.

The threshold representation of a threshold graph may not be unique. 
We prove that it can be chosen such that the clique contains a vertex which is not adjacent to any vertex of the independent set.

\begin{lemma}
    \label{lem:thresholdRepresentation}
    For every threshold graph $G=(V,E)$ there exists a threshold representation $(K,S)$ such that there is a vertex $v\in K$ with $N(v)\cap S=\emptyset$.
\end{lemma}
\begin{proof}
    We fix a threshold representation $(K,S)$ of $G$.
    Suppose for all $v\in K$ holds $N(v)\cap S\neq\emptyset$.
    Then, since $G$ is a threshold graph, there is a vertex $w\in S$ such that $N(w)=K$.
	We obtain a new threshold representation $(K',S')$ of $G$ with $K' \coloneqq K \cup \Set{w}$ and $S' \coloneqq S \setminus \Set{w}$.
    Since $S$ is an independent set, $w$ has no neighbours in $S'$.
\end{proof}

We can now prove that \cref{conj:Tuza} holds for all threshold graphs.

\begin{proof}[of \cref{th:main}]
    Let $G=(K\cup S,E)$ be a threshold graph with $K=\{c_1,\dots,c_k\}$ and $S=\{u_1,\dots,u_s\}$ such that ${N(c_k) \cap S} = \emptyset$. 
    By \cref{lem:thresholdRepresentation}, such a representation exists.
    Let $r \in \Set{1,\dots,s}$ be chosen minimal such that $\Set{c_1,\ldots,c_{\Ceil{ \nicefrac{k}{2}}}} \subseteq N(u_r)$ and let $X$ be the subset $\Set{u_{r}, \dots, u_{s}}$ of $S$  (see \cref{fig:thr-str}).
    Note that $X$ is complete to the set $\Set{c_1,\ldots,c_{\Ceil{\nicefrac{k}{2}}}}$.
    \begin{figure}[htb]
        \centering
        \begin{tikzpicture}[xscale=\ifshort.5\fi\iflong.6\fi,yscale=\ifshort.3\fi\iflong.36\fi]

\def\rndcrn{5pt}
\def\crnshift{8pt}

\tikzstyle{vertex}=[circle,minimum size=4pt,inner sep=0pt,draw,fill]

%%% clique part coordinates

\coordinate (Ksw) at (0,-10.2);
\coordinate (Kne) at (4,0);

\coordinate (Ktopsw) at (0,-5);
\coordinate (Kbotne) at (4,-5);

\coordinate (Klab) at (-1,-5);

\def\Kvx{3.2} % x position of vertices inside K

\coordinate (c1) at (\Kvx,-.5);
\coordinate (c2) at (\Kvx,-1.5);
\coordinate (clhalf) at (\Kvx,-4.5);
\coordinate (cuhalf) at (\Kvx,-5.5);
\coordinate (ck1) at (\Kvx,-8.5);
\coordinate (ck) at (2,-9.5);

%%% stable part coordinates

\coordinate (Ssw) at (6,-8);
\coordinate (Sne) at (10,0);
\coordinate (Snw) at (6,0);

\coordinate (Xnw) at (6,-4);
\coordinate (Xne) at (10,-4);

\coordinate (Xlab) at (9,-6);
\coordinate (Slab) at (11,-4);

\def\Svx{6.8} % x position of vertices inside S

\coordinate (u1) at (\Svx,-.5);
\coordinate (ur) at (\Svx,-3.5);
\coordinate (ur1) at (\Svx,-4.5);
\coordinate (us) at (\Svx,-7.5);

%%% fill K

\draw[thick,fill=black!20,rounded corners=\rndcrn] (Ksw) rectangle (Kne);
%%% draw c_k's neighborhood

\begin{scope}
\draw[fill=black!30] (0,-1)--(0,-9)--(ck)--(4,-9)--(4,-1)--cycle;
\draw[fill=black!30,rounded corners=\rndcrn] (4,-9)--(4,0)--(0,0)--(0,-9);

%\draw[fill=black!30] (4,-9)--(4,0)--(0,0)--(0,-9)--cycle;
\end{scope}

%%% draw full join between X and Ktop

\begin{scope}
\draw[fill=black!20] ([yshift=\crnshift]Ssw)--(Kbotne)--([yshift=-\crnshift]Kne)--(Xnw)--cycle;
\end{scope}

%%% and not so full remaining joins

\begin{scope}
\draw[fill=black!10] (4,-9)--([yshift=\crnshift]Ssw)--(Kbotne)--cycle;
\draw[fill=black!10] ([yshift=-\crnshift]Snw)--(Xnw)--([yshift=-\crnshift]Kne)--cycle;
\end{scope}

%%% draw parts

\draw[thick,rounded corners=\rndcrn] (Ksw) rectangle (Kne);
\draw[thick,rounded corners=\rndcrn] (Ssw) rectangle (Sne);

%%% draw division between Ktop and Kbot and separating X

\draw[dashed] (Ktopsw)--(Kbotne);
\draw[dashed] (Xnw)--(Xne);

%% mark somehow the connections between parts?
%\draw[dashed] (Xnw)--(Kbotne);

%%% label parts

\node at (Klab) {$K$};
\node at (Slab) {$S$};

\node at (Xlab) {$X$};

%%% draw and label vertices

\node[vertex,label=180:$c_1$] at (c1) {};
\node[vertex,label=180:$c_2$] at (c2) {};
\node[vertex,label={[yshift=1pt]180:$c_{\Ceil{\nicefrac{k}{2}}}$}] at (clhalf) {};
\node[vertex,label={[yshift=-2pt]180:$c_{\Ceil{\nicefrac{k}{2}}+1}$}] at (cuhalf) {};
\node[vertex,label=180:$c_{k-1}$] at (ck1) {};
\ifshort
\node[vertex,label={[xshift=1,yshift=-2]180:$c_k$}] at (ck) {};
\fi
\iflong
\node[vertex,label={[xshift=1,yshift=-3]180:$c_k$}] at (ck) {};
\fi

\node at  ($ (c2)!2/5!(clhalf) $) {$\vdots$};
\node at  ($ (ck1)!3/5!(cuhalf) $) {$\vdots$};

\node[vertex,label=0:$u_1$] at (u1) {};
\node[vertex,label=0:$u_{r-1}$] at (ur) {};
\node[vertex,label=0:$u_r$] at (ur1) {};
\node[vertex,label=0:$u_s$] at (us) {};

\node at  ($ (u1)!2/5!(ur) $) {$\vdots$};
\node at  ($ (ur1)!2/5!(us) $) {$\vdots$};
\end{tikzpicture}
        \caption{The structure of threshold graph $G$.}
        \label{fig:thr-str}
    \end{figure}
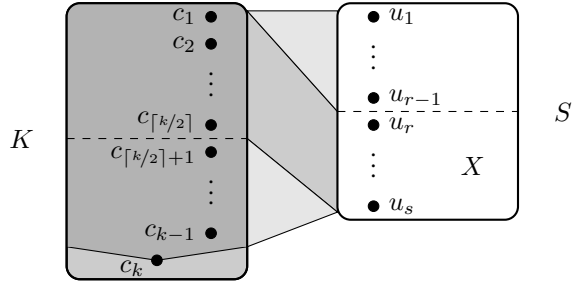
    We distinguish two cases, based on the parity of $k$. 
    First, we focus on the case that $k$ is even.
    In this case we consider two cliques $\Kt$ and $\Kb$ of equal size, induced  by vertices $\Set{c_1,\dots,c_{\nicefrac{k}{2}}}$ and $\Set{c_{\nicefrac{k}{2}+1},\dots,c_k}$, respectively.
    
    We construct a triangle packing $\TPset$ of $G$ using \cref{cor:clique_IS} as follows: we pack the edges of $\Kb$  with vertices in $\Kt$, and the edges of $\Kt$ with vertices in $X$ (see \cref{fig:thr-tp1}).
    
    If $\Abs{X} \geq \frac{k}{2}$, then $\TPset$ is a triangle packing of size $2{\binom{\nicefrac{k}{2}}{2}}$.
    On the other hand, a triangle hitting of size ${\binom{k-1}{2}}$ can be obtained by taking all edges from $K$ except those incident to $c_k$ (see \cref{fig:thr-th1}).
    Thus, we obtain a lower bound on the triangle packing and an upper bound on the triangle hitting yielding:
    \begin{align*}
        \TH{G} \leq \binom{k-1}{2} = \frac{k-2}{2}\cdot(k-1) \leq \frac{k-2}{2}\cdot k = 4{\binom{\nicefrac{k}{2}}{2}} \leq 2\TP{G}.
    \end{align*}
    
\begin{figure}[tb]
     \centering
     \begin{subfigure}[b]{0.49\textwidth}
         \centering
         \begin{tikzpicture}[xscale=\ifshort.5\fi\iflong.6\fi,yscale=\ifshort.3\fi\iflong.36\fi]

\def\rndcrn{5pt}
\def\crnshift{12pt}

\tikzstyle{vertex}=[circle,minimum size=2pt,inner sep=0pt,draw,fill]

%%%%% PACKING

%%% clique part coordinates

\coordinate (Ksw) at (1,-10);
\coordinate (Kne) at (4,0);

\coordinate (Ktopsw) at (1,-5);
\coordinate (Kbotne) at (4,-5);

\coordinate (Klab) at (0,-5);
\coordinate (Ktoplab) at (2.5,-2.5);
\coordinate (Kbotlab) at (2.5,-7.6);

%%% stable part coordinates

\coordinate (Ssw) at (5,-8);
\coordinate (Sne) at (8,0);

\coordinate (Xnw) at (5,-2.5);
\coordinate (Xne) at (8,-2.5);

\coordinate (Xlab) at (6.5,-5.25);
\coordinate (Slab) at (9,-4);

\coordinate (figlab) at (6.5,-9);

%%% fill K

%\draw[thick,fill=black!20] (Ksw) rectangle (Kne);

%%% fill triangle-packings

\begin{scope}[opacity=.5]
\draw[color=blue,fill=blue!50] (Xlab) -- ([shift={(-5pt,-\crnshift-8pt)}] Kne) -- ([shift={(-\rndcrn-5pt,8pt)}] Kbotne) -- cycle;
\end{scope}

\draw[thick,color=blue,fill=blue!20,rounded corners=\rndcrn] ([shift={(5pt,8pt)}] Ktopsw) rectangle ([shift={(-5pt,-8pt)}] Kne);

\begin{scope}[opacity=.5]
\draw[color=blue,fill=blue!50] (Ktoplab) -- ([shift={(\rndcrn+5pt,-8pt)}] Ktopsw) -- ([shift={(-\rndcrn-5pt,-8pt)}] Kbotne) -- cycle;
\end{scope}

\draw[thick,color=blue,fill=blue!20,rounded corners=\rndcrn] ([shift={(5pt,8pt)}] Ksw) rectangle ([shift={(-5pt,-8pt)}] Kbotne);

%%% draw parts

\draw[thick,rounded corners=\rndcrn] (Ksw) rectangle (Kne);
\draw[thick,rounded corners=\rndcrn] (Ssw) rectangle (Sne);

%%% draw division between Ktop and Kbot and separating X

\draw[dashed] (Ktopsw)--(Kbotne);
\draw[dashed] (Xnw)--(Xne);

%%% draw packing from X

\begin{scope}
\def\fromXx{3.6}
\def\fromXy{-0.1}

\coordinate (t1) at (\fromXx,{\fromXy-0.5});
\coordinate (t2) at (\fromXx,{\fromXy-1.25});
\coordinate (t3) at (\fromXx,{\fromXy-2});
\coordinate (t4) at (\fromXx,{\fromXy-2.75});
\coordinate (t5) at (\fromXx,{\fromXy-3.5});
\coordinate (t6) at (\fromXx,{\fromXy-4.25});

\draw[color=blue,fill=blue!20,opacity=0.75] (Xlab) -- (t1) -- (t2)-- (Xlab) -- (t3) -- (t4) -- (Xlab) -- (t5) -- (t6) -- cycle;
\draw[very thick,color=blue] (t1)--(t2);
\draw[very thick,color=blue] (t3)--(t4);
\draw[very thick,color=blue] (t5)--(t6);

\node[vertex] at (t1) {};
\node[vertex] at (t2) {};
\node[vertex] at (t3) {};
\node[vertex] at (t4) {};
\node[vertex] at (t5) {};
\node[vertex] at (t6) {};
\end{scope}

%%% draw packing from Ktop

\begin{scope}
\def\fromKtx{1.125}
\def\fromKty{-5.75}

\coordinate (t1) at ({\fromKtx+.25},\fromKty);
\coordinate (t2) at ({\fromKtx+.7},\fromKty);
\coordinate (t3) at ({\fromKtx+1.15},\fromKty);
\coordinate (t4) at ({\fromKtx+1.6},\fromKty);
\coordinate (t5) at ({\fromKtx+2.05},\fromKty);
\coordinate (t6) at ({\fromKtx+2.5},\fromKty);

\draw[color=blue,fill=blue!20,opacity=0.75] (Ktoplab) -- (t1) -- (t2)-- (Ktoplab) -- (t3) -- (t4) -- (Ktoplab) -- (t5) -- (t6) -- cycle;
\draw[very thick,color=blue] (t1)--(t2);
\draw[very thick,color=blue] (t3)--(t4);
\draw[very thick,color=blue] (t5)--(t6);

\node[vertex] at (t1) {};
\node[vertex] at (t2) {};
\node[vertex] at (t3) {};
\node[vertex] at (t4) {};
\node[vertex] at (t5) {};
\node[vertex] at (t6) {};
\end{scope}

%%% label parts

%\node at (Klab) {$K$};
%\node at (Slab) {$S$};

\node[thick,circle,draw=blue, fill=blue!20, inner sep=0pt,minimum size=25pt] at (Ktoplab) {};
\node[thick,circle,draw=blue, fill=white, inner sep=0pt,minimum size=20pt] at (Xlab) {};

\node at (Kbotlab) {$\Kb$};
\node at (Ktoplab) {$\Kt$};

\node at (Xlab) {$X$};

\node[color=blue] at (figlab) {$\TP{G}$};

\end{tikzpicture}
         \subcaption{triangle packing}
         \label{fig:thr-tp1}
     \end{subfigure}
     \hfill
     \begin{subfigure}[b]{0.49\textwidth}
         \centering
         \begin{tikzpicture}[xscale=\ifshort.5\fi\iflong.6\fi,yscale=\ifshort.3\fi\iflong.36\fi]

\def\rndcrn{5pt}
\def\crnshift{12pt}

\tikzstyle{vertex}=[circle,minimum size=2pt,inner sep=0pt,draw,fill]

%%%%% HITTING

\begin{scope}[xshift=9cm]

%%% clique part coordinates

\coordinate (Ksw) at (1,-10);
\coordinate (Ksw2) at (1,-9);
\coordinate (Kne) at (4,0);

\coordinate (Ktopsw) at (1,-5);
\coordinate (Kbotne) at (4,-5);

\coordinate (Klab) at (2.5,-4.5);
\coordinate (Ktoplab) at (2.5,-2.5);
\coordinate (Kbotlab) at (2.5,-7.6);

\coordinate (ck) at (2.5,-9.5);

%%% stable part coordinates

\coordinate (Ssw) at (5,-8);
\coordinate (Sne) at (8,0);

\coordinate (Xnw) at (5,-2.5);
\coordinate (Xne) at (8,-2.5);

\coordinate (Xlab) at (6.5,-5.25);
\coordinate (Slab) at (9,-4);

\coordinate (figlab) at (6.5,-9);

%%% fill K

%\draw[thick,fill=black!20] (Ksw) rectangle (Kne);

\draw[thick,color=red,fill=red!20,rounded corners=\rndcrn] ([shift={(5pt,8pt)}] Ksw2) rectangle ([shift={(-5pt,-8pt)}] Kne);

%%% draw parts

\draw[thick,rounded corners=\rndcrn] (Ksw) rectangle (Kne);
\draw[thick,rounded corners=\rndcrn] (Ssw) rectangle (Sne);

%%% draw division separating X

\draw[dashed] (Xnw)--(Xne);

%%% label parts

\node at (Klab) {$K$};
%\node at (Slab) {$S$};

%\node at (Kbotlab) {$K_{bot}$};
%\node at (Ktoplab) {$K_{top}$};

\node at (Xlab) {$X$};

\node[color=red] at (figlab) {$\TH{G}$};

%%% draw and label vertices

\node[vertex,label={[yshift=1pt]0:$c_k$}] at (ck) {};

\end{scope}

\end{tikzpicture}
         \subcaption{triangle hitting}
         \label{fig:thr-th1}
     \end{subfigure}
        \caption{The \subref{fig:thr-tp1} triangle packing and \subref{fig:thr-th1} triangle hitting providing the bounds for $|X|\geq\nicefrac{k}{2}$.}
        \label{fig:thr-hp1}
\end{figure}

    If $\Abs{X} < \frac{k}{2}$, then $\TPset$ is of size at least
    \begin{align*}
        \binom{\nicefrac{k}{2}}{2} + \Abs{X} \cdot \Floor{\frac{k}{4}} \geq \binom{\nicefrac{k}{2}}{2} + \Abs{X} \Brace{\frac{k}{4} -\frac{1}{2}}.
    \end{align*}
    On the other hand, the edges inside $\Kt$ and inside $\Kb$ together with all edges between $S$ and $\Kb$ build a triangle hitting of $G$ (cf.\@ \cref{fig:thr-th2}) of size at most
    \begin{align*}
        2\binom{\nicefrac{k}{2}}{2} + \Abs{X}\Brace{\frac{k}{2}- 1}.
    \end{align*}
    Indeed, recall that $c_k$ does not have any neighbours in $S$, therefore we have at most $\Abs{X}\Brace{\frac{k}{2}- 1}$ edges between $X$ and $\Kb$, and by definition of $X$, there are no vertices in $\Kb$ having neighbours in $S \setminus X$.
    Thus, we again obtain a lower bound on the triangle packing and an upper bound on the triangle hitting yielding:
    \begin{align*}
      			\TH{G} \leq 2\binom{\nicefrac{k}{2}}{2} + \Abs{X}\Brace{\frac{k}{2}- 1} = 2 \binom{\nicefrac{k}{2}}{2} + 2\Abs{X}\Brace{\frac{k}{4} -\frac{1}{2}} \leq 2\TP{G}.
      		\end{align*}
\begin{figure}[htb]
     \centering
     \begin{subfigure}[b]{0.49\textwidth}
         \centering
         \begin{tikzpicture}[xscale=\ifshort.5\fi\iflong.6\fi,yscale=\ifshort.3\fi\iflong.36\fi]

\def\rndcrn{5pt}
\def\crnshift{12pt}

\tikzstyle{vertex}=[circle,minimum size=2pt,inner sep=0pt,draw,fill]

%%%%% PACKING

%%% clique part coordinates

\coordinate (Ksw) at (1,-10);
\coordinate (Kne) at (4,0);

\coordinate (Ktopsw) at (1,-5);
\coordinate (Kbotne) at (4,-5);

\coordinate (Klab) at (0,-5);
\coordinate (Ktoplab) at (2.5,-2.5);
\coordinate (Kbotlab) at (2.5,-7.6);

%%% stable part coordinates

\coordinate (Ssw) at (5,-8);
\coordinate (Sne) at (8,0);

\coordinate (Xnw) at (5,-4.5);
\coordinate (Xne) at (8,-4.5);

\coordinate (Xlab) at (6.5,-6.25);
\coordinate (Slab) at (9,-4);

\coordinate (figlab) at (6.5,-9);

%%% fill K

%\draw[thick,fill=black!20] (Ksw) rectangle (Kne);

%%% fill triangle-packings

\begin{scope}[opacity=.5]
\draw[color=blue,fill=blue!50] (Xlab) -- ([shift={(-5pt,-\crnshift-8pt)}] Kne) -- ([shift={(-5pt-\rndcrn,8pt)}] Kbotne) -- cycle;
\end{scope}

\draw[thick,color=blue,fill=blue!20,rounded corners=\rndcrn] ([shift={(15pt,8pt)}] Ktopsw) rectangle ([shift={(-5pt,-8pt)}] Kne);

\begin{scope}[opacity=.5]
\draw[color=blue,fill=blue!50] (Ktoplab) -- ([shift={(\rndcrn+5pt,-8pt)}] Ktopsw) -- ([shift={(-\rndcrn-5pt,-8pt)}] Kbotne) -- cycle;
\end{scope}

\draw[thick,color=blue,fill=blue!20,rounded corners=\rndcrn] ([shift={(5pt,8pt)}] Ksw) rectangle ([shift={(-5pt,-8pt)}] Kbotne);

%%% draw parts

\draw[thick,rounded corners=\rndcrn] (Ksw) rectangle (Kne);
\draw[thick,rounded corners=\rndcrn] (Ssw) rectangle (Sne);

%%% draw division between Ktop and Kbot and separating X

\draw[dashed] (Ktopsw)--(Kbotne);
\draw[dashed] (Xnw)--(Xne);

%%% draw packing from X

\begin{scope}
\def\fromXx{3.6}
\def\fromXy{-0.1}

\coordinate (t1) at (\fromXx,{\fromXy-0.5});
\coordinate (t2) at (\fromXx,{\fromXy-1.25});
\coordinate (t3) at (\fromXx,{\fromXy-2});
\coordinate (t4) at (\fromXx,{\fromXy-2.75});
\coordinate (t5) at (\fromXx,{\fromXy-3.5});
\coordinate (t6) at (\fromXx,{\fromXy-4.25});

\draw[color=blue,fill=blue!20,opacity=0.75] (Xlab) -- (t1) -- (t2)-- (Xlab) -- (t3) -- (t4) -- (Xlab) -- (t5) -- (t6) -- cycle;
\draw[very thick,color=blue] (t1)--(t2);
\draw[very thick,color=blue] (t3)--(t4);
\draw[very thick,color=blue] (t5)--(t6);

\node[vertex] at (t1) {};
\node[vertex] at (t2) {};
\node[vertex] at (t3) {};
\node[vertex] at (t4) {};
\node[vertex] at (t5) {};
\node[vertex] at (t6) {};

\def\fromXx{1.25}

\coordinate (tt1) at (\fromXx,{\fromXy-0.5});
\coordinate (tt2) at (\fromXx,{\fromXy-1.25});
\coordinate (tt3) at (\fromXx,{\fromXy-2});
\coordinate (tt4) at (\fromXx,{\fromXy-2.75});
\coordinate (tt5) at (\fromXx,{\fromXy-3.5});
\coordinate (tt6) at (\fromXx,{\fromXy-4.25});

\draw[very thick,color=orange] (tt1)--(tt2);
\draw[very thick,color=orange] (tt3)--(tt4);
\draw[very thick,color=orange] (tt5)--(tt6);

\node[vertex] at (tt1) {};
\node[vertex] at (tt2) {};
\node[vertex] at (tt3) {};
\node[vertex] at (tt4) {};
\node[vertex] at (tt5) {};
\node[vertex] at (tt6) {};

\end{scope}

%%% draw packing from Ktop

\begin{scope}
\def\fromKtx{1.125}
\def\fromKty{-5.75}

\coordinate (t1) at ({\fromKtx+.25},\fromKty);
\coordinate (t2) at ({\fromKtx+.7},\fromKty);
\coordinate (t3) at ({\fromKtx+1.15},\fromKty);
\coordinate (t4) at ({\fromKtx+1.6},\fromKty);
\coordinate (t5) at ({\fromKtx+2.05},\fromKty);
\coordinate (t6) at ({\fromKtx+2.5},\fromKty);

\draw[color=blue,fill=blue!20,opacity=0.75] (Ktoplab) -- (t1) -- (t2)-- (Ktoplab) -- (t3) -- (t4) -- (Ktoplab) -- (t5) -- (t6) -- cycle;
\draw[very thick,color=blue] (t1)--(t2);
\draw[very thick,color=blue] (t3)--(t4);
\draw[very thick,color=blue] (t5)--(t6);

\node[vertex] at (t1) {};
\node[vertex] at (t2) {};
\node[vertex] at (t3) {};
\node[vertex] at (t4) {};
\node[vertex] at (t5) {};
\node[vertex] at (t6) {};
\end{scope}

%%% label parts

%\node at (Klab) {$K$};
%\node at (Slab) {$S$};

\node[thick,circle,draw=blue, fill=blue!20, inner sep=0pt,minimum size=25pt] at (Ktoplab) {};
\node[thick,circle,draw=blue, fill=white, inner sep=0pt,minimum size=20pt] at (Xlab) {};

\node at (Kbotlab) {$\Kb$};
\node at (Ktoplab) {$\Kt$};

\node at (Xlab) {$X$};

\node[color=blue] at (figlab) {$\TP{G}$};

\end{tikzpicture}
         \subcaption{triangle packing}
         \label{fig:thr-tp2}
     \end{subfigure}
     \hfill
     \begin{subfigure}[b]{0.49\textwidth}
         \centering
         \begin{tikzpicture}[xscale=\ifshort.5\fi\iflong.6\fi,yscale=\ifshort.3\fi\iflong.36\fi]

\def\rndcrn{5pt}
\def\crnshift{12pt}

\tikzstyle{vertex}=[circle,minimum size=2pt,inner sep=0pt,draw,fill]

%%%%% HITTING

\begin{scope}[xshift=9cm]

%%% clique part coordinates

\coordinate (Ksw) at (1,-10);
\coordinate (Ksw2) at (1,-9);
\coordinate (Kne) at (4,0);

\coordinate (Ktopsw) at (1,-5);
\coordinate (Kbotne) at (4,-5);
\coordinate (Kbotse2) at (4,-9);

\coordinate (Klab) at (2.5,-4.5);
\coordinate (Ktoplab) at (2.5,-2.5);
\coordinate (Kbotlab) at (2.5,-7.6);

\coordinate (ck) at (2.5,-9.3);

%%% stable part coordinates

\coordinate (Ssw) at (5,-8);
\coordinate (Sne) at (8,0);

\coordinate (Xnw) at (5,-4.5);
\coordinate (Xne) at (8,-4.5);

\coordinate (Xlab) at (6.5,-6.25);
\coordinate (Slab) at (9,-4);

\coordinate (figlab) at (6.5,-9);

%%% fill K

%\draw[thick,fill=black!20] (Ksw) rectangle (Kne);

\draw[thick,color=red,fill=red!20,rounded corners=\rndcrn] ([shift={(5pt,8pt)}] Ksw) rectangle ([shift={(-5pt,-8pt)}] Kbotne);

\draw[thick,color=red,fill=red!20,rounded corners=\rndcrn] ([shift={(5pt,8pt)}] Ktopsw) rectangle ([shift={(-5pt,-8pt)}] Kne);

%%% draw partial join between X and Kbot

\begin{scope}
\draw[color=red,fill=red!20] ([yshift=\crnshift]Ssw)--(Kbotse2)--(Kbotne)--(Xnw)--cycle;
\end{scope}

%%% draw parts

\draw[thick,rounded corners=\rndcrn] (Ksw) rectangle (Kne);
\draw[thick,rounded corners=\rndcrn] (Ssw) rectangle (Sne);

%%% draw division separating X and parts of K

\draw[dashed] (Ktopsw)--(Kbotne);
\draw[dashed] (Xnw)--(Xne);

%%% label parts

%\node at (Klab) {$K$};
%\node at (Slab) {$S$};

\node at (Kbotlab) {$\Kb$};
\node at (Ktoplab) {$\Kt$};

\node at (Xlab) {$X$};

\node[color=red] at (figlab) {$\TH{G}$};

%%% draw and label vertices

\node[vertex,label={[yshift=1pt]0:$c_k$}] at (ck) {};

\end{scope}

\end{tikzpicture}
         \subcaption{triangle hitting}
         \label{fig:thr-th2}
     \end{subfigure}
        \caption{The \subref{fig:thr-tp2} triangle packing and \subref{fig:thr-th2} triangle hitting providing the bounds when  $|X|<\nicefrac{k}{2}$.}
        \label{fig:thr-hp2}
\end{figure}
    We are left with the case that $k$ is odd.
    We consider the cliques $\Kt$ and $\Kb$ 
    induced by sets $\Set{c_1,\dots,c_{\nicefrac{(k+1)}{2}}}$ and $\Set{c_{\nicefrac{(k+1)}{2}+1},\dots,c_k}$, respectively.

    Again, we look at the size of $X$ and in case it is large, we can derive a similar argument as in the previous case, using \cref{cor:clique_IS}.
    More precisely, assume that $\Abs{X} \geq \frac{k+1}{2}$.
    Then we pack the edges of $\Kb$ into $\binom{\nicefrac{(k-1)}{2}}{2}$ triangles with vertices in $\Kt$, and the edges of $\Kt$ into $\binom{\nicefrac{(k+1)}{2}}{2}$ triangles with vertices in $X$.
    Together, this gives a triangle packing of size
    \begin{align*}
        \binom{\frac{k+1}{2}}{2} +\binom{\frac{k-1}{2}}{2}=\frac{(k-1)^2}{4}.
    \end{align*}
    The triangle hitting 
    again consists of all edges from $K$ except those adjacent to $c_k$, therefore has size $\binom{k-1}{2}$ (recall \cref{fig:thr-hp1}).
    These two bounds together yield:
    \begin{align*}
		\TH{G} \leq \binom{k-1}{2} = \frac{k-1}{2}\cdot\Brace{k-2} \leq \frac{\Brace{k-1}^2}{2} \leq 2\TP{G}.
	\end{align*}
	It remains to consider the case $\Abs{X} < \frac{k+1}{2}$. In order to find a triangle packing, we define $\Kt'$ and $\Kb'$ to be induced by $\Set{c_1,\dots,c_{\nicefrac{(k-1)}{2}}}$ and $\Set{c_{\nicefrac{(k+1)}{2}},\dots,c_k}$, respectively (so $\Kt' = \Kt\setminus\{c_{\nicefrac{(k+1)}{2}}\}$ is of size $\frac{k-1}{2}$ and $\Kb'=\Kb \cup \{c_{\nicefrac{(k+1)}{2}}\}$ is of size $\frac{k+1}{2}$).
    We build a triangle packing analogously to before, using \cref{cor:clique_IS}. 
    The edges of $\Kb'$ can be packed into $\lfloor\frac{\nicefrac{(k+1)}{2}}{2}\rfloor\cdot \frac{k-1}{2}\geq \frac{k-1}{4}\cdot\frac{k-1}{2}$ triangles with vertices in $\Kt'$.
    Moreover, $\min\Set{\Abs{X} \cdot \Floor{\frac{k-1}{4}},\binom{\nicefrac{\Brace{k-1}}{2}}{2}} \geq |X|\frac{k-3}{4}$
    edges of $\Kt'$ can be packed into triangles with vertices in $X$ (see \cref{fig:thr-tp3}).
    This gives a triangle packing of size at least
    \begin{align*}
        \frac{k-1}{2}\cdot\frac{k-1}{4} + \Abs{X}\frac{k-3}{4}.
    \end{align*}
    To find a triangle hitting, we again consider the partition of $K$ into $\Kt$ and $\Kb$.
	We take all edges inside $\Kt$ and inside $\Kb$ together with all edges between $S$ and $\Kb$ (see \cref{fig:thr-th3}).
	Again, recall that $c_k \in \Kb$ does not have any neighbours in $S$, and there are no vertices in $\Kb$ having neighbours in $S \setminus X$.
	Thus, this yields a triangle hitting of size at most.
	\begin{align*}
	    \binom{\frac{k+1}{2}}{2}+\binom{\frac{k-1}{2}}{2}+ \Abs{X} \frac{k-3}{2}.
	\end{align*}
    Therefore, we obtain the following which concludes the proof:
    \begin{align*}
		\TH{G} &\leq \binom{\frac{k+1}{2}}{2}+\binom{\frac{k-1}{2}}{2}+ \Abs{X} \frac{k-3}{2}&\\
           &= \frac{\Brace{k-1}^2}{4}+ \Abs{X} \frac{k-3}{2} = 2\cdot \frac{k-1}{2}\cdot\frac{k-1}{4} + 2\Abs{X}\frac{k-3}{4} \leq 2\TP{G}. \qedhere
	\end{align*}
\end{proof}
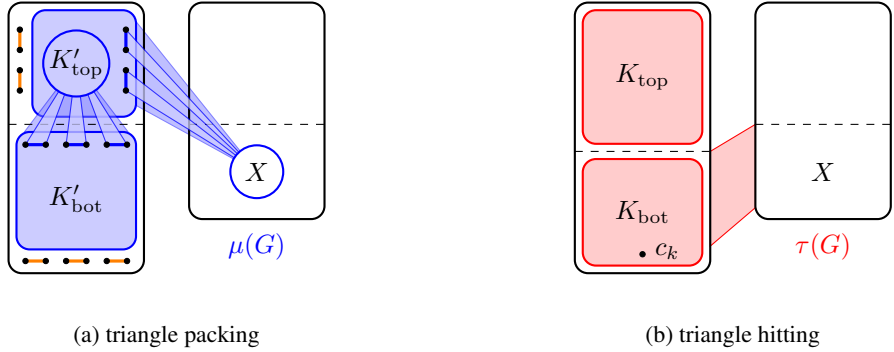
\begin{figure}[ht]
     \centering
     \begin{subfigure}[b]{0.49\textwidth}
         \centering
         \begin{tikzpicture}[xscale=\ifshort.5\fi\iflong.6\fi,yscale=\ifshort.3\fi\iflong.36\fi]

\def\rndcrn{5pt}
\def\crnshift{12pt}

\tikzstyle{vertex}=[circle,minimum size=2pt,inner sep=0pt,draw,fill]

%%%%% PACKING

%%% clique part coordinates

\coordinate (Ksw) at (1,-10);
\coordinate (Kne) at (4,0);

\coordinate (Ktopsw) at (1,-4.5);
\coordinate (Kbotne) at (4,-4.5);
\coordinate (Kbotse2) at (4,-9);

\coordinate (Klab) at (2.5,-4.5);
\coordinate (Ktoplab) at (2.5,-2.25);
\coordinate (Kbotlab) at (2.5,-7.25);

%%% stable part coordinates

\coordinate (Ssw) at (5,-8);
\coordinate (Sne) at (8,0);

\coordinate (Xnw) at (5,-4.5);
\coordinate (Xne) at (8,-4.5);

\coordinate (Xlab) at (6.5,-6.25);
\coordinate (Slab) at (9,-4);

\coordinate (figlab) at (6.5,-9);

%%% fill K

%\draw[thick,fill=black!20] (Ksw) rectangle (Kne);

%%% fill triangle-packings

\begin{scope}[opacity=.5]
\draw[color=blue,fill=blue!50] (Xlab) -- ([shift={(-5pt,-\crnshift-8pt)}] Kne) -- ([shift={(-5pt,\crnshift+8pt)}] Kbotne) -- cycle;
\end{scope}

\draw[thick,color=blue,fill=blue!20,rounded corners=\rndcrn] ([shift={(15pt,8pt)}] Ktopsw) rectangle ([shift={(-5pt,-8pt)}] Kne);

\begin{scope}[opacity=.5]
\draw[color=blue,fill=blue!50] (Ktoplab) -- ([shift={(\rndcrn+5pt,-8pt)}] Ktopsw) -- ([shift={(-\rndcrn-5pt,-8pt)}] Kbotne) -- cycle;
\end{scope}

\draw[thick,color=blue,fill=blue!20,rounded corners=\rndcrn] ([shift={(5pt,25pt)}] Ksw) rectangle ([shift={(-5pt,-8pt)}] Kbotne);

%%% draw parts

\draw[thick,rounded corners=\rndcrn] (Ksw) rectangle (Kne);
\draw[thick,rounded corners=\rndcrn] (Ssw) rectangle (Sne);

%%% draw division between Ktop and Kbot and separating X

\draw[dashed] (Ktopsw)--(Kbotne);
\draw[dashed] (Xnw)--(Xne);

%%% draw packing from X

\begin{scope}
\def\fromXx{3.6}
\def\fromXy{-0.5}

\coordinate (t1) at (\fromXx,{\fromXy-0.5});
\coordinate (t2) at (\fromXx,{\fromXy-1.25});
\coordinate (t3) at (\fromXx,{\fromXy-2});
\coordinate (t4) at (\fromXx,{\fromXy-2.75});
%\coordinate (t5) at (\fromXx,{\fromXy-3.5});
%\coordinate (t6) at (\fromXx,{\fromXy-4.25});

\draw[color=blue,fill=blue!20,opacity=0.75] (Xlab) -- (t1) -- (t2)-- (Xlab) -- (t3) -- (t4) %-- (Xlab) -- (t5) -- (t6) 
-- cycle;
\draw[very thick,color=blue] (t1)--(t2);
\draw[very thick,color=blue] (t3)--(t4);
%\draw[very thick,color=blue] (t5)--(t6);

\node[vertex] at (t1) {};
\node[vertex] at (t2) {};
\node[vertex] at (t3) {};
\node[vertex] at (t4) {};
%\node[vertex] at (t5) {};
%\node[vertex] at (t6) {};

\def\fromXx{1.25}

\coordinate (tt1) at (\fromXx,{\fromXy-0.5});
\coordinate (tt2) at (\fromXx,{\fromXy-1.25});
\coordinate (tt3) at (\fromXx,{\fromXy-2});
\coordinate (tt4) at (\fromXx,{\fromXy-2.75});
%\coordinate (tt5) at (\fromXx,{\fromXy-3.5});
%\coordinate (tt6) at (\fromXx,{\fromXy-4.25});

\draw[very thick,color=orange] (tt1)--(tt2);
\draw[very thick,color=orange] (tt3)--(tt4);
%\draw[very thick,color=orange] (tt5)--(tt6);

\node[vertex] at (tt1) {};
\node[vertex] at (tt2) {};
\node[vertex] at (tt3) {};
\node[vertex] at (tt4) {};
%\node[vertex] at (tt5) {};
%\node[vertex] at (tt6) {};

\end{scope}

%%% draw packing from Ktop

\begin{scope}
\def\fromKtx{1.125}
\def\fromKty{-5.25}

\coordinate (t1) at ({\fromKtx+.25},\fromKty);
\coordinate (t2) at ({\fromKtx+.7},\fromKty);
\coordinate (t3) at ({\fromKtx+1.15},\fromKty);
\coordinate (t4) at ({\fromKtx+1.6},\fromKty);
\coordinate (t5) at ({\fromKtx+2.05},\fromKty);
\coordinate (t6) at ({\fromKtx+2.5},\fromKty);

\draw[color=blue,fill=blue!20,opacity=0.75] (Ktoplab) -- (t1) -- (t2)-- (Ktoplab) -- (t3) -- (t4) -- (Ktoplab) -- (t5) -- (t6) -- cycle;
\draw[very thick,color=blue] (t1)--(t2);
\draw[very thick,color=blue] (t3)--(t4);
\draw[very thick,color=blue] (t5)--(t6);

\node[vertex] at (t1) {};
\node[vertex] at (t2) {};
\node[vertex] at (t3) {};
\node[vertex] at (t4) {};
\node[vertex] at (t5) {};
\node[vertex] at (t6) {};

\def\fromKty{-9.55}

\coordinate (tt1) at ({\fromKtx+.25},\fromKty);
\coordinate (tt2) at ({\fromKtx+.7},\fromKty);
\coordinate (tt3) at ({\fromKtx+1.15},\fromKty);
\coordinate (tt4) at ({\fromKtx+1.6},\fromKty);
\coordinate (tt5) at ({\fromKtx+2.05},\fromKty);
\coordinate (tt6) at ({\fromKtx+2.5},\fromKty);

\draw[very thick,color=orange] (tt1)--(tt2);
\draw[very thick,color=orange] (tt3)--(tt4);
\draw[very thick,color=orange] (tt5)--(tt6);

\node[vertex] at (tt1) {};
\node[vertex] at (tt2) {};
\node[vertex] at (tt3) {};
\node[vertex] at (tt4) {};
\node[vertex] at (tt5) {};
\node[vertex] at (tt6) {};

\end{scope}

%%% label parts

%\node at (Klab) {$K$};
%\node at (Slab) {$S$};

\node[thick,circle,draw=blue, fill=blue!20, inner sep=0pt,minimum size=25pt] at (Ktoplab) {};
\node[thick,circle,draw=blue, fill=white, inner sep=0pt,minimum size=20pt] at (Xlab) {};

\node at (Kbotlab) {$\Kb'$};
\node at (Ktoplab) {$\Kt'$};

\node at (Xlab) {$X$};

\node[color=blue] at (figlab) {$\TP{G}$};

\end{tikzpicture}
         \subcaption{triangle packing}
         \label{fig:thr-tp3}
     \end{subfigure}
     \hfill
     \begin{subfigure}[b]{0.49\textwidth}
         \centering
         \begin{tikzpicture}[xscale=\ifshort.5\fi\iflong.6\fi,yscale=\ifshort.3\fi\iflong.36\fi]

\def\rndcrn{5pt}
\def\crnshift{12pt}

\tikzstyle{vertex}=[circle,minimum size=2pt,inner sep=0pt,draw,fill]

%%%%% HITTING

\begin{scope}[xshift=9cm]

%%% clique part coordinates

\coordinate (Ksw) at (1,-10);
\coordinate (Ksw2) at (1,-9);
\coordinate (Kne) at (4,0);

\coordinate (Ktopsw) at (1,-5.5);
\coordinate (Kbotne) at (4,-5.5);
\coordinate (Kbotse2) at (4,-9);

\coordinate (Klab) at (2.5,-5.5);
\coordinate (Ktoplab) at (2.5,-2.75);
\coordinate (Kbotlab) at (2.5,-7.75);

\coordinate (ck) at (2.5,-9.3);

%%% stable part coordinates

\coordinate (Ssw) at (5,-8);
\coordinate (Sne) at (8,0);

\coordinate (Xnw) at (5,-4.5);
\coordinate (Xne) at (8,-4.5);

\coordinate (Xlab) at (6.5,-6.25);
\coordinate (Slab) at (9,-4);

\coordinate (figlab) at (6.5,-9);

%%% fill K

%\draw[thick,fill=black!20] (Ksw) rectangle (Kne);

\draw[thick,color=red,fill=red!20,rounded corners=\rndcrn] ([shift={(5pt,8pt)}] Ksw) rectangle ([shift={(-5pt,-8pt)}] Kbotne);

\draw[thick,color=red,fill=red!20,rounded corners=\rndcrn] ([shift={(5pt,8pt)}] Ktopsw) rectangle ([shift={(-5pt,-8pt)}] Kne);

%%% draw partial join between X and Kbot

\begin{scope}
\draw[color=red,fill=red!20] ([yshift=\crnshift]Ssw)--(Kbotse2)--(Kbotne)--(Xnw)--cycle;
\end{scope}

%%% draw parts

\draw[thick,rounded corners=\rndcrn] (Ksw) rectangle (Kne);
\draw[thick,rounded corners=\rndcrn] (Ssw) rectangle (Sne);

%%% draw division separating X and parts of K

\draw[dashed] (Ktopsw)--(Kbotne);
\draw[dashed] (Xnw)--(Xne);

%%% label parts

%\node at (Klab) {$K$};
%\node at (Slab) {$S$};

\node at (Kbotlab) {$\Kb$};
\node at (Ktoplab) {$\Kt$};

\node at (Xlab) {$X$};

\node[color=red] at (figlab) {$\TH{G}$};

%%% draw and label vertices

\node[vertex,label={[yshift=1pt]0:$c_k$}] at (ck) {};

\end{scope}

\end{tikzpicture}
         \subcaption{triangle hitting}
         \label{fig:thr-th3}
     \end{subfigure}
        \caption{In \subref{fig:thr-tp3} the triangle packing and in \subref{fig:thr-th2} the triangle hitting providing the bounds for $|K|$ odd and $|X|<\nicefrac{(k+1)}{2}$.}
        \label{fig:thr-hp3}
\end{figure}

\section{Even balanced co-chain graphs}
\begin{figure}[htb]
    \centering
    \begin{tikzpicture}[xscale=\ifshort1\fi\iflong1.2\fi,yscale=\ifshort1\fi\iflong1.2\fi]
	\tikzstyle{vertex}=[shape=circle, fill=black, draw, inner sep=.4mm]
	\tikzstyle{edge}=[draw, gray]
	\pgfdeclarelayer{background}    % declare background layer
\pgfsetlayers{background,main}
	\node (center) at (0,0) {};
	\node (c-center) at ($(center)-(1,0)$) {};
	\node (c2-center) at ($(center)+(1,0)$) {};
	\def\dist{12pt}
		\node (K1label) at ($(c-center)+(-0,{4.5*\dist+4})$) {$K_1$};
		\node (K1label) at ($(c2-center)+(-0,{4.5*\dist+4})$) {$K_2$};
	\foreach \x in {0,...,7}
	{
		\node[vertex] (u-\x) at ($(c-center)+(0.2,{3.5*\dist-\x*\dist})$) {};
		\pgfmathtruncatemacro\y{\x+1}
		\node (u-\x-label) at ($(c-center)+(-0.1,{3.5*\dist-\x*\dist})$) {$c_\y$};
		
		\node[vertex] (d-\x) at ($(c2-center)+(-0.2,{3.5*\dist-\x*\dist})$) {};
		\node (d-\x-label) at ($(c2-center)+(0.15,{3.5*\dist-\x*\dist})$) {$d_\y$};
	}
% newly added
%    \foreach \start/\end in {0/3,4/7}{
%	    \draw [black,decorate,decoration={brace,amplitude=5pt, raise=3.5*\dist},xshift=0,yshift=0pt] (u-\end) -- (u-\start) node [black,midway,xshift=-4.5*\dist] {$\ell$};
%	     \draw [black,decorate,decoration={brace,amplitude=5pt, raise=3.5*\dist},xshift=0,yshift=0pt] (d-\start) -- (d-\end) node [black,midway,xshift=4.5*\dist] {$\ell$};
%	}
	
%
	\draw[rounded corners] ($(c-center)-(\dist,4*\dist)$) rectangle ($(c-center)+(\dist,4*\dist)$);
	\draw[rounded corners] ($(c2-center)-(\dist,4*\dist)$) rectangle ($(c2-center)+(\dist,4*\dist)$);
	%draw ($(c-center)-(\dist,\dist)$) rectangle ($(c-center)+(\dist,4*\dist)$);
	
	\begin{pgfonlayer}{background}
		\begin{scope}
			\clip ($(c-center)-(\dist,\dist)$) rectangle ($(c-center)+(\dist,4*\dist)$);
			\fill[green!40!black,opacity=0.3,rounded corners] ($(c-center)-(\dist,4*\dist)$) rectangle ($(c-center)+(\dist,4*\dist)$);
		\end{scope}
		\node (X_1) at ($(c-center)+(-1.9*\dist,0.3*\dist)$) {\textcolor{green!40!black}{$X_1$}}; 
	\end{pgfonlayer}
	
	\begin{pgfonlayer}{background}
		\begin{scope}
			\clip ($(c2-center)-(\dist,4*\dist)$) rectangle ($(c2-center)+(\dist,\dist)$);
			\fill[green!40!black,opacity=0.3, rounded corners] ($(c2-center)-(\dist,4*\dist)$) rectangle ($(c2-center)+(\dist,4*\dist)$);
		\end{scope}
		\node (X_2) at ($(c2-center)+(1.9*\dist,-0.3*\dist)$) {\textcolor{green!40!black}{$X_2$}}; 
	\end{pgfonlayer}

	\draw[dashed] ($(c-center)-(\dist,0)$) to ($(c-center)+(\dist,0)$);
	\draw[dashed] ($(c2-center)-(\dist,0)$) to ($(c2-center)+(\dist,0)$);
	
	\node at ($(c-center)+(-2.2*\dist,2*\dist)$) {$\Kto$};
	\node at ($(c-center)+(-2.2*\dist,-2*\dist)$) {$\Kbo$};
	
	\node at ($(c2-center)+(2.2*\dist,2*\dist)$) {$\Ktt$};
	\node at ($(c2-center)+(2.2*\dist,-2*\dist)$) {$\Kbt$};
	
	%%%% EDGES BETWEEN CLIQUES 
    \foreach \x in {0,...,3}{
	      \draw[edge] (u-0) -- (d-\x);
	      \draw[edge] (u-\x) -- (d-3);
	   \foreach \y in {4,...,7}{
	       \draw[edge] (u-\x) -- (d-\y);
	   }   
	}
	\draw[edge] (u-4) -- (d-4);
	\draw[edge] (u-4) -- (d-5);
	\draw[edge] (u-4) -- (d-6);
	\draw[edge] (u-4) -- (d-7);
	\draw[edge] (u-1) -- (d-2);
	\draw[edge] (u-2) -- (d-2);
	\draw[edge] (u-5) -- (d-5);
	\draw[edge] (u-5) -- (d-6);
	\draw[edge] (u-5) -- (d-7);
	\draw[edge] (u-6) -- (d-7);

\end{tikzpicture}
    \caption{An example of an even balanced co-chain graph with $\ell=4$ (omitting the edges inside the cliques $K_1$ and $K_2$).
    }
    \label{fig:co-chain-pic}
\end{figure}
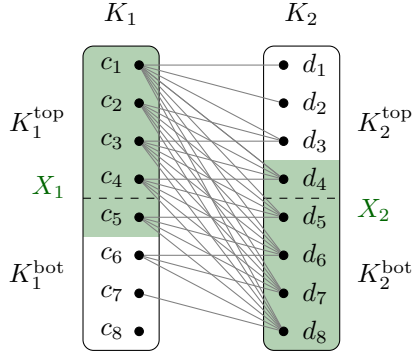

In this section we prove \cref{thm:cochain}.
To this end let $G$ be an even balanced co-chain graph and $\Brace{K_1,K_2}$ its \emph{co-chain representation}.
Recall that $K_1$ and $K_2$ are of same size which is divisible by 4, for the rest of the section let $\Abs{K_1}=\Abs{K_2}=2\ell$ for $\ell$ even. We identify $K_1$ and $K_2$ with the cliques $G[K_1]$ and $G[K_2]$. 
See \cref{fig:co-chain-pic} for an illustration.

We prove that Tuza's conjecture holds for this graph class.

\begin{proof}[of \cref{thm:cochain}]
Note that in the case $\ell=2$ we get an $8$-vertex graph which is either a clique, or  has average degree less than $7$, so
this case is covered by \cite{puleo2015tuza}. Therefore in the following we assume that $\ell\geq 4$. 
  
Similarly to threshold graphs, we use $\Kto,\Kbo$ for the top and the bottom half of $K_1$, respectively, and similarly $\Ktt,\Kbt$ for the top and the bottom half of $K_2$.
Let $X_1\subseteq K_1,X_2\subseteq K_2$ be the sets defined as follows: $c\in X_1$ if $\Kbt \subseteq N[c]$, and $d\in X_2$ if $\Kto \subseteq N[d]$.
See~\cref{fig:co-chain-pic} for an illustration.
We denote $x_1=|X_1|$ and $x_2=|X_2|$.
By definition, $x_1\ge\ell$ implies that the set $X_1\supseteq \Kto$ is complete to $\Kbt$. Consequently, $x_2\ge\ell$. Similarly, $x_2\ge\ell$ implies $x_1\ge\ell$.
Therefore, $x_1\ge\ell$ if and only if $x_2\ge\ell$.
We assume without loss of generality throughout the entire proof that $x_1 \geq x_2$.
We split the analysis into two main cases.

\subsection{The case \texorpdfstring{$x_1,x_2\le\ell$}{x1,x2<=l}}
In this case $X_1\subseteq\Kto$ and $X_2\subseteq \Kbt$.
Suppose there is an edge $cd$ with $c \in K_1\setminus X_1$ and $d \in \Ktt$, then $c$ is adjacent to all the vertices in $\Kbt$ and so $c\in X_1$, which yields a contradiction.
Similarly, there are no edges between $\Kbo$ and $K_2\setminus X_2$.
In particular, there are no edges between $\Ktt$ and $\Kbo$.

We choose a triangle hitting $\THset$ obtained by taking all edges within $\Kto$, $\Ktt$, $\Kbo$, and $\Kbt$, as well as edges between $X_1$ and $\Kbt$, and between $X_2$ and $\Kto$ as illustrated in \cref{fig:cobip-hitting-1}.
Observe now that in the graph $G-\THset$ vertices in $X_1$ only have neighbours in the independent set $\Kbo \cup \Ktt$, vertices in $\Kto\setminus X_1$ only have neighbours in the independent set $\Kbo \cup \Kbt\setminus X_2$, while vertices in $\Kbo$ only have neighbours in the independent set $\Kto \cup X_2$. Therefore the set $\THset$ is indeed a triangle hitting of $G$.
	\begin{figure}[htb]
        \centering
        \begin{tikzpicture}[xscale=\ifshort.5\fi\iflong.6\fi,yscale=\ifshort.3\fi\iflong.36\fi]

\def\rndcrn{5pt}
\def\crnshift{12pt}

\tikzstyle{vertex}=[circle,minimum size=2pt,inner sep=0pt,draw,fill]

%%%%% HITTING

%%% clique 1 part coordinates

\coordinate (K1sw) at (1,-10);
\coordinate (K1ne) at (4,0);

\coordinate (K1topsw) at (1,-5);
\coordinate (K1botne) at (4,-5);

\coordinate (K1lab) at (0,-5);
\coordinate (K1toplab) at (2.5,-2.3);
\coordinate (K1botlab) at (2.5,-7.6);

\coordinate (KX1nw) at (1,-3);
\coordinate (KX1ne) at (4,-3);
\coordinate (KX1lab) at (2.5,-4);

%%% clique 2 part coordinates

\coordinate (K2sw) at (5,-10);
\coordinate (K2ne) at (8,0);

\coordinate (K2topsw) at (5,-5);
\coordinate (K2botne) at (8,-5);

\coordinate (K2lab) at (9,-5);
\coordinate (K2toplab) at (6.5,-2.3);
\coordinate (K2botlab) at (6.5,-7.6);

\coordinate (X2nw) at (5,-7);
\coordinate (X2ne) at (8,-7);
\coordinate (X2lab) at (6.5,-8.5);

\coordinate (figlab) at (4.5,-11);

%%% draw division between Ktop and Kbot and separating X

\draw[dashed] (K1topsw)--(K1botne);
\draw[dashed] (K2topsw)--(K2botne);

\draw[color=green!40!black] (KX1nw)--(KX1ne);
\draw[color=green!40!black] (X2nw)--(X2ne);

%fill hitting

\draw[thick,color=red,fill=red!20,rounded corners=\rndcrn] ([shift={(5pt,8pt)}] K1sw) rectangle ([shift={(-5pt,-8pt)}] K1botne);

\draw[thick,color=red,fill=red!40,opacity=0.5,rounded corners=\rndcrn] ([shift={(5pt,8pt)}] K2sw) rectangle ([shift={(-5pt,-8pt)}] K2botne);
\draw[thick,color=red,rounded corners=\rndcrn] ([shift={(5pt,8pt)}] K2sw) rectangle ([shift={(-5pt,-8pt)}] K2botne);

\draw[thick,color=red,fill=red!40,opacity=0.5,rounded corners=\rndcrn] ([shift={(5pt,8pt)}] K1topsw) rectangle ([shift={(-5pt,-8pt)}] K1ne);
\draw[thick,color=red,rounded corners=\rndcrn] ([shift={(5pt,8pt)}] K1topsw) rectangle ([shift={(-5pt,-8pt)}] K1ne);

\draw[thick,color=red,fill=red!20,rounded corners=\rndcrn] ([shift={(5pt,8pt)}] K2topsw) rectangle ([shift={(-5pt,-8pt)}] K2ne);

%%% draw partial join between X and Kbot

\begin{scope}
\draw[color=red,fill=red!40] ([yshift=\crnshift]K2sw)--(KX1ne)--([yshift=-\crnshift]K1ne)--(X2nw)--cycle;
\end{scope}

\begin{scope}
\draw[color=red,fill=red!20] ([yshift=\crnshift]K2sw)--(KX1ne)--(K1botne)--cycle;
\end{scope}

\begin{scope}
\draw[color=red,fill=red!20] (K2topsw)--([yshift=-\crnshift]K1ne)--(X2nw)--cycle;
\end{scope}

%%% draw parts

\draw[thick,rounded corners=\rndcrn] (K1sw) rectangle (K1ne);
\draw[thick,rounded corners=\rndcrn] (K2sw) rectangle (K2ne);

%%% label parts

\node at (K1botlab) {$\Kbo$};
\node at (K1toplab) {$\Kto$};

\node at (K2botlab) {$\Kbt$};
\node at (K2toplab) {$\Ktt$};

\end{tikzpicture}
	    \caption{The triangle hitting used in the case $x_1,x_2\leq\ell$.}
        \label{fig:cobip-hitting-1}
    \end{figure}
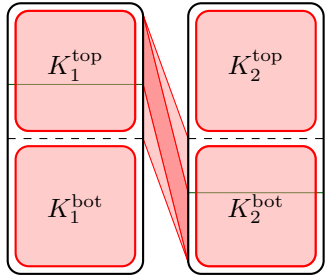

Therefore,
\[
  \TH{G} \le |\THset| = 4 \binom{\ell}{2}+\ell x_1+\ell x_2-x_1 x_2=4 \binom{\ell}{2}+\ell x_1+(\ell-x_1) x_2.
\]
Indeed, we note that we counted edges between $X_1$ and $X_2$ once in term $\ell x_1$ and once in term $\ell x_2$ which we compensate by subtracting the last term $x_1 x_2$.

Let us now create a sufficiently large triangle packing.
First, we pack all edges of $\Kbo$ with vertices in $\Kto$ and also all edges of $\Ktt$ with vertices in $\Kbt$; we denote the set of these triangles by $A$ (see \cref{fig:cob-tpA}).
By \cref{lem:clique_IS}, $A$ contains $2\binom{\ell}{2}$ triangles.
Observe that $2|A|-|\THset| = -\ell x_1 - (\ell-x_1)x_2$.
First, we sort out the single case where $x_1=\ell$, and, in consequence, $x_2=\ell$ by definition of $X_1$ and $X_2$ together with the assumption that $x_2\le\ell$.

\subsubsection{The subcase \texorpdfstring{$x_1=x_2=\ell$}{x1=x2=l}}

In this case, $|\THset|=4\binom{\ell}{2}+\ell^2$.
As $\Kto\cup\Kbt$ is a clique, by \cref{lem:clique} we can pack at least $\frac{1}{3}\Brace{\binom{2\ell}{2}-\ell-1}$ triangles in it.
Together with $A$, we obtain a triangle packing $\TPset$. 
If $\ell\geq 5$, then $2\TPset-\THset\ge\frac{2}{3}\left(\binom{2\ell}{2}-\ell-1\right)-\ell^2=\frac{1}{3}\left(\ell^2-4\ell-2\right)\ge 0$.
If $\ell=4$, \cref{lem:clique} gives us a stronger bound without the term $-2$, leading to $2\TPset-\THset\ge\frac{1}{3}\left(\ell^2-4\ell\right)= 0$.
Both cases imply $2\TP{G} \geq \TH{G}$.

\begin{figure}[htb]
     \centering
     \begin{subfigure}[b]{0.49\textwidth}
         \centering
         \begin{tikzpicture}[xscale=\ifshort.5\fi\iflong.6\fi,yscale=\ifshort.3\fi\iflong.36\fi]

\def\rndcrn{5pt}
\def\crnshift{12pt}

\tikzstyle{vertex}=[circle,minimum size=2pt,inner sep=0pt,draw,fill]

%%%%% PACKING

%%% clique 1 part coordinates

\coordinate (K1sw) at (1,-10);
\coordinate (K1ne) at (4,0);

\coordinate (K1topsw) at (1,-5);
\coordinate (K1botne) at (4,-5);

\coordinate (K1lab) at (0,-5);
\coordinate (K1toplab) at (2.5,-2.5);
\coordinate (K1botlab) at (2.5,-7.6);

%%% clique 2 part coordinates

\coordinate (K2sw) at (5,-10);
\coordinate (K2ne) at (8,0);

\coordinate (K2topsw) at (5,-5);
\coordinate (K2botne) at (8,-5);

\coordinate (K2lab) at (9,-5);
\coordinate (K2toplab) at (6.5,-2.5);
\coordinate (K2botlab) at (6.5,-7.6);

\coordinate (figlab) at (4.5,-11);

%%% fill triangle-packings

\begin{scope}[opacity=.5]
\draw[color=blue,fill=blue!50] (K1toplab) -- ([shift={(\rndcrn+5pt,-8pt)}] K1topsw) -- ([shift={(-\rndcrn-5pt,-8pt)}] K1botne) -- cycle;
\end{scope}

\draw[thick,color=blue,fill=blue!20,rounded corners=\rndcrn] ([shift={(5pt,8pt)}] K1sw) rectangle ([shift={(-5pt,-8pt)}] K1botne);

\begin{scope}[opacity=.5]
\draw[color=blue,fill=blue!50] (K2botlab) -- ([shift={(\rndcrn+5pt,+8pt)}] K2topsw) -- ([shift={(-\rndcrn-5pt,+8pt)}] K2botne) -- cycle;
\end{scope}

\draw[thick,color=blue,fill=blue!20,rounded corners=\rndcrn] ([shift={(5pt,8pt)}] K2topsw) rectangle ([shift={(-5pt,-8pt)}] K2ne);

%%% draw parts

\draw[thick,rounded corners=\rndcrn] (K1sw) rectangle (K1ne);
\draw[thick,rounded corners=\rndcrn] (K2sw) rectangle (K2ne);

%%% draw division between Ktop and Kbot and separating X

\draw[dashed] (K1topsw)--(K1botne);
\draw[dashed] (K2topsw)--(K2botne);

%%% draw packing from K1top

\begin{scope}
\def\fromKtx{1.125}
\def\fromKty{-5.75}

\coordinate (t1) at ({\fromKtx+.25},\fromKty);
\coordinate (t2) at ({\fromKtx+.7},\fromKty);
\coordinate (t3) at ({\fromKtx+1.15},\fromKty);
\coordinate (t4) at ({\fromKtx+1.6},\fromKty);
\coordinate (t5) at ({\fromKtx+2.05},\fromKty);
\coordinate (t6) at ({\fromKtx+2.5},\fromKty);

\draw[color=blue,fill=blue!20,opacity=0.75] (K1toplab) -- (t1) -- (t2)-- (K1toplab) -- (t3) -- (t4) -- (K1toplab) -- (t5) -- (t6) -- cycle;
\draw[very thick,color=blue] (t1)--(t2);
\draw[very thick,color=blue] (t3)--(t4);
\draw[very thick,color=blue] (t5)--(t6);

\node[vertex] at (t1) {};
\node[vertex] at (t2) {};
\node[vertex] at (t3) {};
\node[vertex] at (t4) {};
\node[vertex] at (t5) {};
\node[vertex] at (t6) {};
\end{scope}

%%% draw packing from K2bot

\begin{scope}
\def\fromKtx{5.125}
\def\fromKty{-4.25}

\coordinate (t1) at ({\fromKtx+.25},\fromKty);
\coordinate (t2) at ({\fromKtx+.7},\fromKty);
\coordinate (t3) at ({\fromKtx+1.15},\fromKty);
\coordinate (t4) at ({\fromKtx+1.6},\fromKty);
\coordinate (t5) at ({\fromKtx+2.05},\fromKty);
\coordinate (t6) at ({\fromKtx+2.5},\fromKty);

\draw[color=blue,fill=blue!20,opacity=0.75] (K2botlab) -- (t1) -- (t2)-- (K2botlab) -- (t3) -- (t4) -- (K2botlab) -- (t5) -- (t6) -- cycle;
\draw[very thick,color=blue] (t1)--(t2);
\draw[very thick,color=blue] (t3)--(t4);
\draw[very thick,color=blue] (t5)--(t6);

\node[vertex] at (t1) {};
\node[vertex] at (t2) {};
\node[vertex] at (t3) {};
\node[vertex] at (t4) {};
\node[vertex] at (t5) {};
\node[vertex] at (t6) {};
\end{scope}

%%% label parts

\node[thick,circle,draw=blue, fill=white, inner sep=0pt,minimum size=25pt] at (K1toplab) {};

\node at (K1botlab) {$\Kbo$};
\node at (K1toplab) {$\Kto$};

\node[thick,circle,draw=blue, fill=white, inner sep=0pt,minimum size=25pt] at (K2botlab) {};

\node at (K2botlab) {$\Kbt$};
\node at (K2toplab) {$\Ktt$};

%\node[color=blue] at (figlab) {$A$};

\end{tikzpicture}
         \subcaption{set $A$ of triangles}
         \label{fig:cob-tpA}
     \end{subfigure}
     \hfill
     \begin{subfigure}[b]{0.49\textwidth}
         \centering
         \begin{tikzpicture}[xscale=\ifshort.5\fi\iflong.6\fi,yscale=\ifshort.3\fi\iflong.36\fi]

\def\rndcrn{5pt}
\def\crnshift{12pt}

\tikzstyle{vertex}=[circle,minimum size=2pt,inner sep=0pt,draw,fill]

%%%%% PACKING

%%% clique 1 part coordinates

\coordinate (K1sw) at (1,-10);
\coordinate (K1ne) at (4,0);

\coordinate (K1topsw) at (1,-5);
\coordinate (K1botne) at (4,-5);

\coordinate (K1lab) at (0,-5);
\coordinate (K1toplab) at (2.5,-2.5);
\coordinate (K1botlab) at (2.5,-7.6);

\coordinate (X1sw) at (1,-3);
\coordinate (X1se) at (4,-3);
\coordinate (X1lab) at (2.5,-1.5);

%%% clique 2 part coordinates

\coordinate (K2sw) at (5,-10);
\coordinate (K2ne) at (8,0);

\coordinate (K2topsw) at (5,-5);
\coordinate (K2botne) at (8,-5);

\coordinate (K2lab) at (9,-5);
\coordinate (K2toplab) at (6.5,-2.5);
\coordinate (K2botlab) at (6.5,-7.6);

\coordinate (figlab) at (4.5,-11);

%%% fill triangle-packings

\begin{scope}[opacity=.5]
\draw[color=violet,fill=violet!50] (X1lab) -- ([shift={(15pt,-\crnshift-8pt)}] K2topsw) -- ([shift={(5pt,\crnshift+8pt)}] K2sw) -- cycle;
\end{scope}

\draw[thick,color=violet,fill=violet!20,rounded corners=\rndcrn] ([shift={(5pt,8pt)}] K2sw) rectangle ([shift={(-15pt,-8pt)}] K2botne);

%%% draw parts

\draw[thick,rounded corners=\rndcrn] (K1sw) rectangle (K1ne);
\draw[thick,rounded corners=\rndcrn] (K2sw) rectangle (K2ne);

%%% draw division between Ktop and Kbot and separating X

\draw[dashed] (K1topsw)--(K1botne);
\draw[dashed] (K2topsw)--(K2botne);

\draw[color=green!40!black] (X1sw)--(X1se);

%%% draw packing from X

\begin{scope}
\def\fromXx{5.4}
\def\fromXy{-5.4}

\coordinate (t1) at (\fromXx,{\fromXy-0.25});
\coordinate (t2) at (\fromXx,{\fromXy-1});
\coordinate (t3) at (\fromXx,{\fromXy-1.75});
\coordinate (t4) at (\fromXx,{\fromXy-2.5});
\coordinate (t5) at (\fromXx,{\fromXy-3.25});
\coordinate (t6) at (\fromXx,{\fromXy-4});

\draw[color=violet,fill=violet!20,opacity=0.75] (X1lab) -- (t1) -- (t2)-- (X1lab) -- (t3) -- (t4) -- (X1lab) -- (t5) -- (t6) 
-- cycle;
\draw[very thick,color=violet] (t1)--(t2);
\draw[very thick,color=violet] (t3)--(t4);
\draw[very thick,color=violet] (t5)--(t6);

\node[vertex] at (t1) {};
\node[vertex] at (t2) {};
\node[vertex] at (t3) {};
\node[vertex] at (t4) {};
\node[vertex] at (t5) {};
\node[vertex] at (t6) {};

\def\fromXx{7.75}

\coordinate (tt1) at (\fromXx,{\fromXy-0.25});
\coordinate (tt2) at (\fromXx,{\fromXy-1});
\coordinate (tt3) at (\fromXx,{\fromXy-1.75});
\coordinate (tt4) at (\fromXx,{\fromXy-2.5});
\coordinate (tt5) at (\fromXx,{\fromXy-3.25});
\coordinate (tt6) at (\fromXx,{\fromXy-4});

\draw[very thick,color=orange] (tt1)--(tt2);
\draw[very thick,color=orange] (tt3)--(tt4);
\draw[very thick,color=orange] (tt5)--(tt6);

\node[vertex] at (tt1) {};
\node[vertex] at (tt2) {};
\node[vertex] at (tt3) {};
\node[vertex] at (tt4) {};
\node[vertex] at (tt5) {};
\node[vertex] at (tt6) {};

\end{scope}

%%% label parts

\node[thick,circle,draw=violet, fill=white, inner sep=0pt,minimum size=20pt] at (X1lab) {};

\node at (K2botlab) {$\Kbt$};

\node at (X1lab) {$X_1$};

%\node[color=violet] at (figlab) {$B$};

\end{tikzpicture}
         \subcaption{set $B$ of triangles}
         \label{fig:cob-tpB}
     \end{subfigure}
     \\[.5cm]
     \begin{subfigure}[b]{0.49\textwidth}
         \centering
         \begin{tikzpicture}[xscale=\ifshort.5\fi\iflong.6\fi,yscale=\ifshort.3\fi\iflong.36\fi]

\def\rndcrn{5pt}
\def\crnshift{12pt}

\tikzstyle{vertex}=[circle,minimum size=2pt,inner sep=0pt,draw,fill]

%%%%% PACKING

%%% clique 1 part coordinates

\coordinate (K1sw) at (1,-10);
\coordinate (K1ne) at (4,0);

\coordinate (K1topsw) at (1,-5);
\coordinate (K1botne) at (4,-5);

\coordinate (K1lab) at (0,-5);
\coordinate (K1toplab) at (2.5,-2.5);
\coordinate (K1botlab) at (2.5,-7.6);

\coordinate (X1sw) at (1,-3);
\coordinate (X1se) at (4,-3);
\coordinate (X1lab) at (2.5,-1.7);

\coordinate (KX1nw) at (1,-3);
\coordinate (KX1ne) at (4,-3);
\coordinate (KX1lab) at (2.5,-4);

%%% clique 2 part coordinates

\coordinate (K2sw) at (5,-10);
\coordinate (K2ne) at (8,0);

\coordinate (K2topsw) at (5,-5);
\coordinate (K2botne) at (8,-5);

\coordinate (K2lab) at (9,-5);
\coordinate (K2toplab) at (6.5,-2.5);
\coordinate (K2botlab) at (6.5,-7.6);

\coordinate (figlab) at (4.5,-11);

%%% fill triangle-packings

\begin{scope}[opacity=.5]
\draw[color=olive,fill=olive!50] (KX1lab) -- ([shift={(\rndcrn+5pt,+8pt)}] X1sw) -- ([shift={(-\rndcrn-5pt,+8pt)}] X1se) -- cycle;
\end{scope}

\draw[thick,color=olive,fill=olive!20,rounded corners=\rndcrn] ([shift={(5pt,5pt)}] X1sw) rectangle ([shift={(-5pt,-23pt)}] K1ne);

%%% draw parts

\draw[thick,rounded corners=\rndcrn] (K1sw) rectangle (K1ne);
\draw[thick,rounded corners=\rndcrn] (K2sw) rectangle (K2ne);

%%% draw division between Ktop and Kbot and separating X

\draw[dashed] (K1topsw)--(K1botne);
\draw[dashed] (K2topsw)--(K2botne);

\draw[color=green!40!black] (X1sw)--(X1se);

%%% draw packing from K1-X1

\begin{scope}
\def\fromKtx{1.125}
\def\fromKty{-2.5}

\coordinate (t1) at ({\fromKtx+.25},\fromKty);
\coordinate (t2) at ({\fromKtx+.7},\fromKty);
\coordinate (t3) at ({\fromKtx+1.15},\fromKty);
\coordinate (t4) at ({\fromKtx+1.6},\fromKty);
\coordinate (t5) at ({\fromKtx+2.05},\fromKty);
\coordinate (t6) at ({\fromKtx+2.5},\fromKty);

\draw[color=olive,fill=olive!20,opacity=0.75] (KX1lab) -- (t1) -- (t2)-- (KX1lab) -- (t3) -- (t4) -- (KX1lab) -- (t5) -- (t6) -- cycle;
\draw[very thick,color=olive] (t1)--(t2);
\draw[very thick,color=olive] (t3)--(t4);
\draw[very thick,color=olive] (t5)--(t6);

\node[vertex] at (t1) {};
\node[vertex] at (t2) {};
\node[vertex] at (t3) {};
\node[vertex] at (t4) {};
\node[vertex] at (t5) {};
\node[vertex] at (t6) {};

\def\fromKty{-0.4}

\coordinate (t1) at ({\fromKtx+.25},\fromKty);
\coordinate (t2) at ({\fromKtx+.7},\fromKty);
\coordinate (t3) at ({\fromKtx+1.15},\fromKty);
\coordinate (t4) at ({\fromKtx+1.6},\fromKty);
\coordinate (t5) at ({\fromKtx+2.05},\fromKty);
\coordinate (t6) at ({\fromKtx+2.5},\fromKty);

\draw[very thick,color=orange] (t1)--(t2);
\draw[very thick,color=orange] (t3)--(t4);
\draw[very thick,color=orange] (t5)--(t6);

\node[vertex] at (t1) {};
\node[vertex] at (t2) {};
\node[vertex] at (t3) {};
\node[vertex] at (t4) {};
\node[vertex] at (t5) {};
\node[vertex] at (t6) {};
\end{scope}

%%% label parts

\node[thick,rectangle,rounded corners=\rndcrn,draw=olive, fill=white, inner sep=0pt,minimum width=30pt, minimum height=12pt] at (KX1lab) {};

\node at (KX1lab) {\tiny $\Kto\!\!\setminus\!\! X_1$};
\node at (X1lab) {$X_1$};

%\node[color=olive] at (figlab) {$C$};

\end{tikzpicture}
         \subcaption{set $C$ of triangles}
         \label{fig:cob-tpC}
     \end{subfigure}
     \hfill
     \begin{subfigure}[b]{0.49\textwidth}
         \centering
         \begin{tikzpicture}[xscale=\ifshort.5\fi\iflong.6\fi,yscale=\ifshort.3\fi\iflong.36\fi]

\def\rndcrn{5pt}
\def\crnshift{12pt}

\tikzstyle{vertex}=[circle,minimum size=2pt,inner sep=0pt,draw,fill]

%%%%% PACKING

%%% clique 1 part coordinates

\coordinate (K1sw) at (1,-10);
\coordinate (K1ne) at (4,0);

\coordinate (K1topsw) at (1,-5);
\coordinate (K1botne) at (4,-5);

\coordinate (K1lab) at (0,-5);
\coordinate (K1toplab) at (2.5,-2.5);
\coordinate (K1botlab) at (2.5,-7.6);

\coordinate (KX1nw) at (1,-3);
\coordinate (KX1ne) at (4,-3);
\coordinate (KX1lab) at (2.5,-4);

%%% clique 2 part coordinates

\coordinate (K2sw) at (5,-10);
\coordinate (K2ne) at (8,0);

\coordinate (K2topsw) at (5,-5);
\coordinate (K2botne) at (8,-5);

\coordinate (K2lab) at (9,-5);
\coordinate (K2toplab) at (6.5,-2.5);
\coordinate (K2botlab) at (6.5,-7.6);

\coordinate (X2nw) at (5,-7);
\coordinate (X2ne) at (8,-7);
\coordinate (X2lab) at (6.5,-8.5);

\coordinate (figlab) at (4.5,-11);

%%% fill triangle-packings

\begin{scope}[opacity=.5]
\draw[color=cyan,fill=cyan!50] (X2lab) -- ([shift={(-5pt,-\crnshift-8pt)}] KX1ne) -- ([shift={(-20pt,\crnshift+8pt)}] K1botne) -- cycle;
\end{scope}

\draw[thick,color=cyan,fill=cyan!20,rounded corners=\rndcrn] ([shift={(15pt,8pt)}] K1topsw) rectangle ([shift={(-5pt,-8pt)}] KX1ne);

%%% draw parts

\draw[thick,rounded corners=\rndcrn] (K1sw) rectangle (K1ne);
\draw[thick,rounded corners=\rndcrn] (K2sw) rectangle (K2ne);

%%% draw division between Ktop and Kbot and separating X

\draw[dashed] (K1topsw)--(K1botne);
\draw[dashed] (K2topsw)--(K2botne);

\draw[color=green!40!black] (KX1nw)--(KX1ne);
\draw[color=green!40!black] (X2nw)--(X2ne);

%%% draw packing from X

\begin{scope}
\def\fromXx{3.6}
\def\fromXy{-3.4}

\coordinate (t1) at (\fromXx,{\fromXy-0.25});
\coordinate (t2) at (\fromXx,{\fromXy-1});
%\coordinate (t3) at (\fromXx,{\fromXy-1.75});
%\coordinate (t4) at (\fromXx,{\fromXy-2.5});
%\coordinate (t5) at (\fromXx,{\fromXy-3.25});
%\coordinate (t6) at (\fromXx,{\fromXy-4});

\draw[color=cyan,fill=cyan!20,opacity=0.75] (X2lab) -- (t1) -- (t2)-- (X2lab) %-- (t3) -- (t4) -- (X1lab) -- (t5) -- (t6) 
-- cycle;
\draw[very thick,color=cyan] (t1)--(t2);
%\draw[very thick,color=violet] (t3)--(t4);
%\draw[very thick,color=violet] (t5)--(t6);

\node[vertex] at (t1) {};
\node[vertex] at (t2) {};
%\node[vertex] at (t3) {};
%\node[vertex] at (t4) {};
%\node[vertex] at (t5) {};
%\node[vertex] at (t6) {};

\def\fromXx{1.25}

\coordinate (tt1) at (\fromXx,{\fromXy-0.25});
\coordinate (tt2) at (\fromXx,{\fromXy-1});
%\coordinate (tt3) at (\fromXx,{\fromXy-1.75});
%\coordinate (tt4) at (\fromXx,{\fromXy-2.5});
%\coordinate (tt5) at (\fromXx,{\fromXy-3.25});
%\coordinate (tt6) at (\fromXx,{\fromXy-4});

\draw[very thick,color=orange] (tt1)--(tt2);
%\draw[very thick,color=orange] (tt3)--(tt4);
%\draw[very thick,color=orange] (tt5)--(tt6);

\node[vertex] at (tt1) {};
\node[vertex] at (tt2) {};
%\node[vertex] at (tt3) {};
%\node[vertex] at (tt4) {};
%\node[vertex] at (tt5) {};
%\node[vertex] at (tt6) {};

\end{scope}

%%% label parts

\node[thick,circle,draw=cyan, fill=white, inner sep=0pt,minimum size=20pt] at (X2lab) {};

\node at (KX1lab) {\tiny $\Kto\!\!\setminus\!\! X_1$};
\node at (X2lab) {$X_2$};

%\node[color=cyan] at (figlab) {$D$};

\end{tikzpicture}
         \subcaption{set $D$ of triangles}
         \label{fig:cob-tpD}
     \end{subfigure}
        \caption{Triangles in \subref{fig:cob-tpA} $A$, \subref{fig:cob-tpB} $B$, \subref{fig:cob-tpC} $C$, and \subref{fig:cob-tpD} $D$ in the case $x_1,x_2\leq \ell$.}
        \label{fig:cob-tp}
\end{figure}
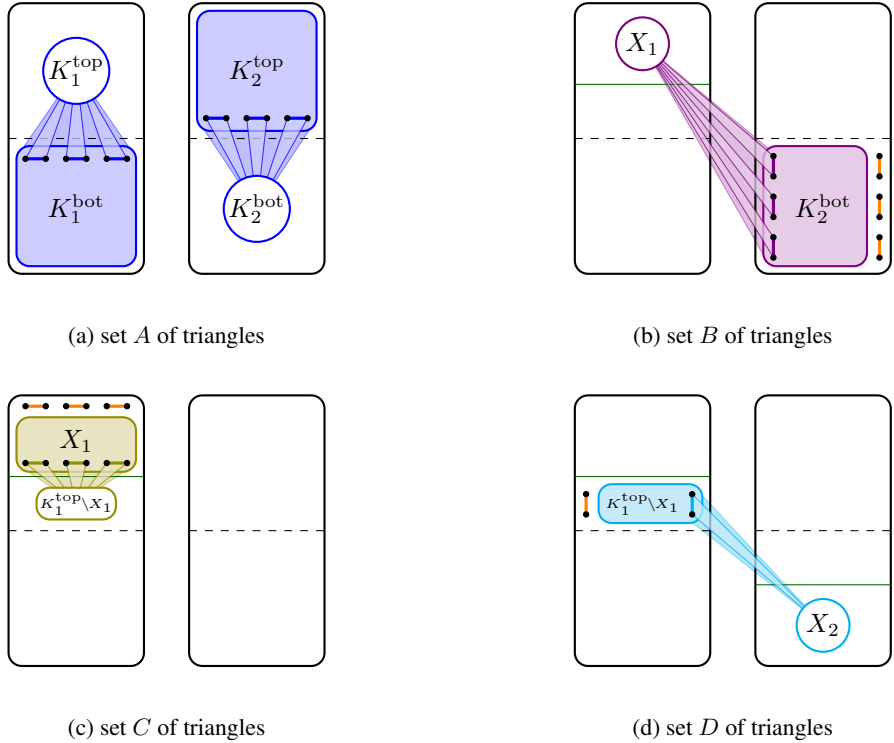

\subsubsection{The subcase \texorpdfstring{$x_1,x_2<\ell$}{x1,x2<l}}
Now, we consider the remaining case where $x_1<\ell$, and, in consequence, $x_2<\ell$.

We choose a triangle packing $\TPset$ as follows (see \cref{fig:cob-tp}).
We take the set $A$ of triangles as defined before.
Recall that $2|A|-|\THset|=-\ell x_1 - (\ell-x_1)x_2$.
We create a set $B$ of triangles by packing edges of $\Kbt$ with vertices in $X_1$.
By \cref{cor:clique_IS}(\ref{cor:clique_IS:between}) and as $x_1<\ell$, $B$ is of size $\ell/2\cdot x_1$.
We create another set of triangles $C$ by packing edges of $X_1$ with vertices of $\Kto\setminus X_1$.
Next, let $D$ be the set of triangles created by packing edges of $\Kto\setminus X_1$ with vertices in $X_2$. 
It is clear that all triangles in $\TPset = A \cup B \cup C \cup D$ are mutually edge-disjoint, therefore $\TPset$ is indeed a triangle packing.

%%%%%x_1 even
Let us first settle the case that $x_1$ is even.
As $2(|A|+|B|)-|\THset|=-(\ell-x_1)x_2$ if $x_1<\ell$, it remains to show that $2\Abs{\TPset \setminus \Brace{A \cup B}}  = 2(\Abs{C}+\Abs{D}) \geq {{(\ell-x_1)x_2}}$.

If $\ell-x_1> x_2$, then $2|D|=(\ell-x_1)x_2$ by~\cref{cor:clique_IS}(\ref{cor:clique_IS:between}).
So, assume that $\ell-x_1\le x_2$.
Consequently, $\ell-x_1\le x_1$ and thus ${\ell}/{2}\leq x_1$.
If $x_1 = \ell/2$, then, by $x_1\ge x_2\ge \ell/2$, we have $x_2=\ell/2$ as well.
Thus, as $\ell \geq 4$, $2(|C|+|D|)-(\ell-x_1)x_2=4\binom{\ell/2}{2}-\ell^2/4=\ell(\ell-4)/4\geq 0$.
For $\ell-x_1<x_1$ we get $2|C|=x_1(\ell-x_1)\ge x_2(\ell-x_1)$.
Therefore, we always have $2\Abs{C \cup D} \geq {{(\ell-x_1)x_2}}$ for even $x_1$, and so $2\TP{G}\ge2\TPset\ge\THset\ge\TH{G}$.

%%%%%%x_1 odd
In case $x_1$ is odd, we add one additional triangle to our triangle packing as follows.
Note that if there is no edge between $\Kbo$ and $\Kbt$, then all edges between $\Kto$ and $\Ktt$ hit all triangles between $K_1$ and $K_2$, therefore taking these edges instead of edges between $\Kto$ and $\Kbt$ creates a triangle hitting $\THset'$ of size at most $4\binom{\ell}{2}+x_1\ell$ as all the edges between $\Kto$ and $\Ktt$ have one endpoint in $X_1$.
As $x_1<\ell$, we obtain $2\TP{G} \geq 2(|A|+|B|)\geq|\THset'|\geq\TH{G}$.
Thus we can assume that there is at least one edge $uv$ with $u\in \Kbo$ and $v\in \Kbt$.

Note in particular that $v\in X_2$ as every edge between $\Kbo$ and $\Kbt$ has one endpoint in $X_2$.
Observe that $|\Kto\setminus X_1|=\ell-x_1$ is odd, so there exists an unpacked matching between $\Kto\setminus X_1$ and $X_2$ (not containing edges used in triangles from set $D$).
Indeed, each maximal matching in $\Kto\setminus X_1$ constructed according to \cref{thm:mat_decomp} omits a different vertex $u_1\in\Kto\setminus X_1$, so after the matching is fully joined with a vertex $u_2\in X_2$, as in \cref{cor:clique_IS}, the edge $u_1u_2$ remains unpacked. A collection of all such edges gives the desired matching.
Let $w\in \Kto\setminus X_1$ be a vertex such that $wv$ is an edge of the mentioned unpacked matching. 
Finally, as $\ell$ is even, a star with center in $\Kto$ is not used in any triangle in $A$, by \cref{lem:clique_IS}. 
Note that the center of this star can be chosen arbitrarily among vertices of $\Kto$ by \cref{lem:clique_IS}; let us choose $w$ to be the center.
Therefore, $uvw$ is a triangle which is edge-disjoint with every triangle in $A\cup B\cup C\cup D$ and we may set $\TPset^{\operatorname{odd}}=\TPset \cup\{uvw\}$ for odd $x_1$.

Recall that $2(|A|+|B|)-|\THset|=-(\ell-x_1)x_2$. Similarly as before, we need to prove that \[2\Abs{\TPset^{\operatorname{odd}} \setminus \Brace{A \cup B}} = 2(\Abs{C} + \Abs{D} + 1) \geq {{(\ell-x_1)x_2}}.\]

\begin{description}
\item If $\ell-x_1\le x_2$, then again $\ell-x_1\le x_1$ and thus ${\ell}/{2}\leq x_1$.
    The case $\ell/2=x_1$ can be handled exactly as in the even case.
    So assume further $\ell-x_1<x_1$, then using \cref{cor:clique_IS} we obtain
    $2(|C|+|D|)=(x_1-1)(\ell-x_1)+2\binom{\ell-x_1}{2}=(x_1-1)(\ell-x_1)+(\ell-x_1)(\ell-x_1-1)=(\ell-x_1)(\ell-2)$.
    Consequently, $2(|C|+|D|+1) - (\ell-x_1)x_2 = 2+(\ell-x_1)(\ell-2-x_2)$.
    Observe that, for $x_2\le\ell-2$, we already get $(\ell-x_1)(\ell-2-x_2)\ge0$.
    We have $x_1 = \ell - 1$ because $x_1$ is odd and $\ell$ is even.
    For $x_2=\ell-1$, we have $x_1=\ell-1$ because $x_2\le x_1 < \ell$.
    Thus $2+(\ell-x_1)(\ell-2-x_2) = 2+ 1\cdot(-1)\geq 0$.
    Therefore, we obtain $2\Brace{\Abs{C} + \Abs{D} + 1} \geq {(\ell-x_1)x_2}$.
    \newpage
\item If $\ell-x_1> x_2$,
    then $2|D|=(\ell-x_1-1)x_2=(\ell-x_1)x_2-x_2$.
    Hence in this case, $D$ alone does not suffice as it is missing $x_2$ triangles.
    We therefore need $2\Abs{C}+2 \geq x_2$.
    We use \cref{cor:clique_IS} to analyse the size of $C$.
    
    If $x_1\le\ell-x_1$, then $2|C|+2-x_2\ge x_1(x_1-1)-x_2+2\ge (x_2-1)^2+1\ge 1$ as $x_1(x_1-1)\ge x_2(x_2-1)$.
    If $x_1>\ell-x_1$, then, 
    $2|C|+2-x_2=(x_1-1)(\ell-x_1)-x_2+2\ge x_1-x_2+1\geq 1$, as $\ell-x_1\geq 1$ and $x_1\ge x_2$.
    So in both cases we obtain $2\Abs{C}+2 \geq x_2 +1 \geq x_2$.
\end{description}
We conclude that $2\TP{G}\ge2\TPset^{\operatorname{odd}}\ge\THset\ge\TH{G}$.

\subsection{The case \texorpdfstring{$x_1>\ell$}{x1>l} and \texorpdfstring{$x_2\ge\ell$}{x2≥l}}
	\begin{figure}[ht]
        \centering
        \begin{tikzpicture}[xscale=\ifshort.5\fi\iflong.6\fi,yscale=\ifshort.3\fi\iflong.36\fi]

\def\rndcrn{5pt}
\def\crnshift{12pt}

\tikzstyle{vertex}=[circle,minimum size=2pt,inner sep=0pt,draw,fill]

%%%%% HITTING

%%% clique 1 part coordinates

\coordinate (K1sw) at (1,-10);
\coordinate (K1ne) at (4,0);

\coordinate (K1topsw) at (1,-5);
\coordinate (K1botne) at (4,-5);

\coordinate (K1lab) at (0,-5);
\coordinate (K1toplab) at (2.5,-2.3);
\coordinate (K1botlab) at (2.5,-7.6);

\coordinate (KX1nw) at (1,-7);
\coordinate (KX1ne) at (4,-7);
\coordinate (KX1lab) at (2.5,-3.5);

%%% clique 2 part coordinates

\coordinate (K2sw) at (5,-10);
\coordinate (K2ne) at (8,0);

\coordinate (K2topsw) at (5,-5);
\coordinate (K2botne) at (8,-5);

\coordinate (K2lab) at (9,-5);
\coordinate (K2toplab) at (6.5,-2.3);
\coordinate (K2botlab) at (6.5,-7.6);

\coordinate (X2nw) at (5,-3);
\coordinate (X2ne) at (8,-3);
\coordinate (X2lab) at (6.5,-6.5);

\coordinate (figlab) at (4.5,-11);

%%% draw division between Ktop and Kbot and separating X

\draw[dashed] (K1topsw)--(K1botne);
\draw[dashed] (K2topsw)--(K2botne);

\draw[color=green!40!black] (KX1nw)--(KX1ne);
\draw[color=green!40!black] (X2nw)--(X2ne);

%fill hitting

\draw[thick,color=red,fill=red!20,rounded corners=\rndcrn] ([shift={(5pt,8pt)}] K2sw) rectangle ([shift={(-5pt,-8pt)}] K2botne);

\draw[thick,color=red,fill=red!40,opacity=0.5,rounded corners=\rndcrn] ([shift={(5pt,8pt)}] K1sw) rectangle ([shift={(-5pt,-8pt)}] K1botne);
\draw[thick,color=red,rounded corners=\rndcrn] ([shift={(5pt,8pt)}] K1sw) rectangle ([shift={(-5pt,-8pt)}] K1botne);

\draw[thick,color=red,fill=red!40,opacity=0.5,rounded corners=\rndcrn] ([shift={(5pt,8pt)}] K2topsw) rectangle ([shift={(-5pt,-8pt)}] K2ne);
\draw[thick,color=red,rounded corners=\rndcrn] ([shift={(5pt,8pt)}] K2topsw) rectangle ([shift={(-5pt,-8pt)}] K2ne);

\draw[thick,color=red,fill=red!20,rounded corners=\rndcrn] ([shift={(5pt,8pt)}] K1topsw) rectangle ([shift={(-5pt,-8pt)}] K1ne);

%%% draw partial join between X and Kbot

\begin{scope}

%\draw[color=red,fill=red!40,opacity=0.5]
\draw[color=red,fill=red!40]([yshift=\crnshift]K2sw)--(K1botne)--([yshift=-\crnshift]K1ne)--(K2topsw)--cycle;
\draw[color=red] ([yshift=\crnshift]K2sw)--(K1botne)--([yshift=-\crnshift]K1ne)--(K2topsw)--cycle;
\end{scope}

\begin{scope}
\draw[color=red,fill=red!40,opacity=0.5] (KX1ne)--(K1botne)--(X2nw)--(K2topsw)--cycle;
\draw[color=red] (KX1ne)--(K1botne)--(X2nw)--(K2topsw)--cycle;
\end{scope}

%%% draw parts

\draw[thick,rounded corners=\rndcrn] (K1sw) rectangle (K1ne);
\draw[thick,rounded corners=\rndcrn] (K2sw) rectangle (K2ne);

%%% label parts

\node at (K1botlab) {$\Kbo$};
\node at (K1toplab) {$\Kto$};

\node at (K2botlab) {$\Kbt$};
\node at (K2toplab) {$\Ktt$};

\end{tikzpicture}
	    \caption{The triangle hitting used in the case $x_1>\ell$ and $x_2\ge \ell$.}
        \label{fig:cobip-hitting-2}
    \end{figure}

We choose a triangle hitting $\THset$ obtained by taking all edges within $\Kto$, $\Kbo$, $\Ktt$ and $\Kbt$ as well as edges between $\Kto$ and $\Kbt$ and between $\Kbo$ and $\Ktt$ (cf.~\cref{fig:cobip-hitting-2}). The graph $G - \THset$ is bipartite, thus $\THset$ is indeed a triangle hitting in $G$. We have
\begin{align*}
    \Abs{\THset} &= 4 \binom{\ell}{2} + \ell^2 + \Abs{\Fkt{E}{\Ktt,\Kbo}}
    \leq 3\ell^2 - 2\ell + (x_1-\ell)(x_2-\ell).
\end{align*}

\begin{figure}[htb]
     \centering
     \begin{subfigure}[b]{0.49\textwidth}
         \centering
         \begin{tikzpicture}[xscale=\ifshort.5\fi\iflong.6\fi,yscale=\ifshort.3\fi\iflong.36\fi]

\def\rndcrn{5pt}
\def\crnshift{12pt}

\tikzstyle{vertex}=[circle,minimum size=2pt,inner sep=0pt,draw,fill]

%%%%% PACKING

%%% clique 1 part coordinates

\coordinate (K1sw) at (1,-10);
\coordinate (K1ne) at (4,0);

\coordinate (K1topsw) at (1,-5);
\coordinate (K1botne) at (4,-5);

\coordinate (K1lab) at (0,-5);
\coordinate (K1toplab) at (2.5,-2.5);
\coordinate (K1botlab) at (2.5,-7.6);

%%% clique 2 part coordinates

\coordinate (K2sw) at (5,-10);
\coordinate (K2ne) at (8,0);

\coordinate (K2topsw) at (5,-5);
\coordinate (K2botne) at (8,-5);

\coordinate (K2lab) at (9,-5);
\coordinate (K2toplab) at (6.5,-2.5);
\coordinate (K2botlab) at (6.5,-7.6);

\coordinate (figlab) at (4.5,-11);

%%% fill triangle-packings

\begin{scope}[opacity=.5]
\draw[color=blue,fill=blue!50] (K1toplab) -- ([shift={(\rndcrn+5pt,-8pt)}] K1topsw) -- ([shift={(-\rndcrn-5pt,-8pt)}] K1botne) -- cycle;
\end{scope}

\draw[thick,color=blue,fill=blue!20,rounded corners=\rndcrn] ([shift={(5pt,8pt)}] K1sw) rectangle ([shift={(-5pt,-8pt)}] K1botne);

\begin{scope}[opacity=.5]
\draw[color=blue,fill=blue!50] (K2botlab) -- ([shift={(\rndcrn+5pt,+8pt)}] K2topsw) -- ([shift={(-\rndcrn-5pt,+8pt)}] K2botne) -- cycle;
\end{scope}

\draw[thick,color=blue,fill=blue!20,rounded corners=\rndcrn] ([shift={(5pt,8pt)}] K2topsw) rectangle ([shift={(-5pt,-8pt)}] K2ne);

\begin{scope}[opacity=.5]
\draw[color=blue,fill=blue!50] (K2botlab) -- ([shift={(-5pt,-\crnshift-26pt)}] K1ne) -- ([shift={(-8pt,\crnshift+8pt)}] K1botne) -- cycle;
\end{scope}

\draw[thick,color=blue,fill=blue!20,rounded corners=\rndcrn] ([shift={(5pt,8pt)}] K1topsw) rectangle ([shift={(-5pt,-8pt)}] K1ne);

%%% draw parts

\draw[thick,rounded corners=\rndcrn] (K1sw) rectangle (K1ne);
\draw[thick,rounded corners=\rndcrn] (K2sw) rectangle (K2ne);

%%% draw division between Ktop and Kbot and separating X

\draw[dashed] (K1topsw)--(K1botne);
\draw[dashed] (K2topsw)--(K2botne);

%%% draw packing from K1top

\begin{scope}
\def\fromKtx{1.125}
\def\fromKty{-5.75}

\coordinate (t1) at ({\fromKtx+.25},\fromKty);
\coordinate (t2) at ({\fromKtx+.7},\fromKty);
\coordinate (t3) at ({\fromKtx+1.15},\fromKty);
\coordinate (t4) at ({\fromKtx+1.6},\fromKty);
\coordinate (t5) at ({\fromKtx+2.05},\fromKty);
\coordinate (t6) at ({\fromKtx+2.5},\fromKty);

\draw[color=blue,fill=blue!20,opacity=0.75] (K1toplab) -- (t1) -- (t2)-- (K1toplab) -- (t3) -- (t4) -- (K1toplab) -- (t5) -- (t6) -- cycle;
\draw[very thick,color=blue] (t1)--(t2);
\draw[very thick,color=blue] (t3)--(t4);
\draw[very thick,color=blue] (t5)--(t6);

\node[vertex] at (t1) {};
\node[vertex] at (t2) {};
\node[vertex] at (t3) {};
\node[vertex] at (t4) {};
\node[vertex] at (t5) {};
\node[vertex] at (t6) {};
\end{scope}

%%% draw packing from K2bot

\begin{scope}
\def\fromKtx{5.125}
\def\fromKty{-4.25}

\coordinate (t1) at ({\fromKtx+.25},\fromKty);
\coordinate (t2) at ({\fromKtx+.7},\fromKty);
\coordinate (t3) at ({\fromKtx+1.15},\fromKty);
\coordinate (t4) at ({\fromKtx+1.6},\fromKty);
\coordinate (t5) at ({\fromKtx+2.05},\fromKty);
\coordinate (t6) at ({\fromKtx+2.5},\fromKty);

\draw[color=blue,fill=blue!20,opacity=0.75] (K2botlab) -- (t1) -- (t2)-- (K2botlab) -- (t3) -- (t4) -- (K2botlab) -- (t5) -- (t6) -- cycle;
\draw[very thick,color=blue] (t1)--(t2);
\draw[very thick,color=blue] (t3)--(t4);
\draw[very thick,color=blue] (t5)--(t6);

\node[vertex] at (t1) {};
\node[vertex] at (t2) {};
\node[vertex] at (t3) {};
\node[vertex] at (t4) {};
\node[vertex] at (t5) {};
\node[vertex] at (t6) {};
\end{scope}

%%% draw another packing from K2bot

\begin{scope}
\def\fromXx{3.6}
\def\fromXy{-1.1}

\coordinate (t1) at (\fromXx,{\fromXy-0.25});
\coordinate (t2) at (\fromXx,{\fromXy-1});
\coordinate (t3) at (\fromXx,{\fromXy-1.75});
\coordinate (t4) at (\fromXx,{\fromXy-2.5});
%\coordinate (t5) at (\fromXx,{\fromXy-3.25});
%\coordinate (t6) at (\fromXx,{\fromXy-4});

\draw[color=blue,fill=blue!20,opacity=0.75] (K2botlab) -- (t1) -- (t2)-- (K2botlab) -- (t3) -- (t4) %-- (K2botlab) -- (t5) -- (t6) 
-- cycle;
\draw[very thick,color=blue] (t1)--(t2);
\draw[very thick,color=blue] (t3)--(t4);
%\draw[very thick,color=blue] (t5)--(t6);

\node[vertex] at (t1) {};
\node[vertex] at (t2) {};
\node[vertex] at (t3) {};
\node[vertex] at (t4) {};
%\node[vertex] at (t5) {};
%\node[vertex] at (t6) {};

\end{scope}

%%% label parts

\node[thick,circle,draw=blue, fill=blue!20, inner sep=0pt,minimum size=25pt] at (K1toplab) {};

\node at (K1botlab) {$\Kbo$};
\node at (K1toplab) {$\Kto$};

\node[thick,circle,draw=blue, fill=white, inner sep=0pt,minimum size=25pt] at (K2botlab) {};

\node at (K2botlab) {$\Kbt$};
\node at (K2toplab) {$\Ktt$};

\node[vertex, xshift=-17pt] at  (K2botlab) {};
\node[ xshift=-17pt,yshift=-4pt] at  (K2botlab) {$v$};

%\node[color=blue] at (figlab) {$A$};

\end{tikzpicture}
         \subcaption{set $A'$ of triangles}
         \label{fig:cob-tp2-1}
     \end{subfigure}
     \hfill
     \begin{subfigure}[b]{0.49\textwidth}
         \centering
         \begin{tikzpicture}[xscale=\ifshort.5\fi\iflong.6\fi,yscale=\ifshort.3\fi\iflong.36\fi]

\def\rndcrn{5pt}
\def\crnshift{12pt}

\tikzstyle{vertex}=[circle,minimum size=2pt,inner sep=0pt,draw,fill]

%%%%% PACKING

%%% clique 1 part coordinates

\coordinate (K1sw) at (1,-10);
\coordinate (K1ne) at (4,0);

\coordinate (K1topsw) at (1,-5);
\coordinate (K1botne) at (4,-5);

\coordinate (K1lab) at (0,-5);
\coordinate (K1toplab) at (2.5,-2.5);
\coordinate (K1botlab) at (2.5,-7.6);

%\coordinate (X1sw) at (1,-6);
%\coordinate (X1se) at (4,-6);
%\coordinate (KX1lab) at (2.5,-5);
\coordinate (X1sw) at (1,-8);
\coordinate (X1se) at (4,-8);
\coordinate (KX1lab) at (2.5,-6.5);

%%% clique 2 part coordinates

\coordinate (K2sw) at (5,-10);
\coordinate (K2ne) at (8,0);

\coordinate (K2topsw) at (5,-5);
\coordinate (K2botne) at (8,-5);

\coordinate (K2lab) at (9,-5);
\coordinate (K2toplab) at (6.5,-2.5);
\coordinate (K2botlab) at (6.5,-7.6);

\coordinate (figlab) at (4.5,-11);

%%% fill triangle-packings

\begin{scope}[opacity=.5]
\draw[color=violet,fill=violet!50] (KX1lab) -- ([shift={(5pt,-\crnshift-3pt)}] K2topsw) -- ([shift={(5pt,\crnshift+5pt)}] K2sw) -- cycle;
\end{scope}

\draw[thick,color=violet,fill=violet!20,rounded corners=\rndcrn] ([shift={(5pt,8pt)}] K2sw) rectangle ([shift={(-15pt,-8pt)}] K2botne);

%%% draw parts

\draw[thick,rounded corners=\rndcrn] (K1sw) rectangle (K1ne);
\draw[thick,rounded corners=\rndcrn] (K2sw) rectangle (K2ne);

%%% draw division between Ktop and Kbot and separating X

\draw[dashed] (K1topsw)--(K1botne);
\draw[dashed] (K2topsw)--(K2botne);

\draw[color=green!40!black] (X1sw)--(X1se);

%%% draw packing from X

\begin{scope}
\def\fromXx{5.4}
\def\fromXy{-5.4}

\coordinate (t1) at (\fromXx,{\fromXy-0.25});
\coordinate (t2) at (\fromXx,{\fromXy-1});
\coordinate (t3) at (\fromXx,{\fromXy-1.75});
\coordinate (t4) at (\fromXx,{\fromXy-2.5});
\coordinate (t5) at (\fromXx,{\fromXy-3.25});
\coordinate (t6) at (\fromXx,{\fromXy-4});

\draw[color=violet,fill=violet!20,opacity=0.75] (KX1lab) -- (t1) -- (t2)-- (KX1lab) -- (t3) -- (t4) -- (KX1lab) -- (t5) -- (t6) 
-- cycle;
\draw[very thick,color=violet] (t1)--(t2);
\draw[very thick,color=violet] (t3)--(t4);
\draw[very thick,color=violet] (t5)--(t6);

\node[vertex] at (t1) {};
\node[vertex] at (t2) {};
\node[vertex] at (t3) {};
\node[vertex] at (t4) {};
\node[vertex] at (t5) {};
\node[vertex] at (t6) {};

\def\fromXx{7.75}

\coordinate (tt1) at (\fromXx,{\fromXy-0.25});
\coordinate (tt2) at (\fromXx,{\fromXy-1});
\coordinate (tt3) at (\fromXx,{\fromXy-1.75});
\coordinate (tt4) at (\fromXx,{\fromXy-2.5});
\coordinate (tt5) at (\fromXx,{\fromXy-3.25});
\coordinate (tt6) at (\fromXx,{\fromXy-4});

\draw[very thick,color=orange] (tt1)--(tt2);
\draw[very thick,color=orange] (tt3)--(tt4);
\draw[very thick,color=orange] (tt5)--(tt6);

\node[vertex] at (tt1) {};
\node[vertex] at (tt2) {};
\node[vertex] at (tt3) {};
\node[vertex] at (tt4) {};
\node[vertex] at (tt5) {};
\node[vertex] at (tt6) {};

\end{scope}

%%% label parts

%\node[thick,circle,draw=violet, fill=white, inner sep=0pt,minimum size=20pt] at (X1lab) {};

%\node[thick,rectangle,rounded corners=\rndcrn,draw=violet, fill=white, inner sep=0pt,minimum width=35pt, minimum height=12pt] at (KX1lab) {};
\node[thick,rectangle,rounded corners=\rndcrn,draw=violet, fill=white, inner sep=0pt,minimum width=37pt, minimum height=20pt] at (KX1lab) {};

\node at (K2botlab) {$\Kbt$};

\node at (KX1lab) {\tiny $X_1\cap \Kbo$};

%\node[color=violet] at (figlab) {$B$};

\end{tikzpicture}
         \subcaption{set $B'$ of triangles}
         \label{fig:cob-tp2-2}
     \end{subfigure}
     \\[.5cm]
     \begin{subfigure}[b]{0.49\textwidth}
         \centering
         \begin{tikzpicture}[xscale=\ifshort.5\fi\iflong.6\fi,yscale=\ifshort.3\fi\iflong.36\fi]

\def\rndcrn{5pt}
\def\crnshift{12pt}

\tikzstyle{vertex}=[circle,minimum size=2pt,inner sep=0pt,draw,fill]

%%%%% PACKING

%%% clique 1 part coordinates

\coordinate (K1sw) at (1,-10);
\coordinate (K1ne) at (4,0);

\coordinate (K1topsw) at (1,-5);
\coordinate (K1botne) at (4,-5);
\coordinate (K1botsw) at (1,-5);

\coordinate (K1lab) at (0,-5);
\coordinate (K1toplab) at (2.5,-2.5);
\coordinate (K1botlab) at (2.5,-7.6);

%%% clique 2 part coordinates

\coordinate (K2sw) at (5,-10);
\coordinate (K2ne) at (8,0);

\coordinate (K2topsw) at (5,-5);
\coordinate (K2botne) at (8,-5);

\coordinate (K2lab) at (9,-5);
\coordinate (K2toplab) at (6.5,-2.5);
\coordinate (K2botlab) at (6.5,-7.6);

\coordinate (figlab) at (4.5,-11);

\coordinate (X2sw) at (5,-2.3);
\coordinate (X2se) at (8,-2.3);
\coordinate (KX2lab) at (6.5,-3.6);

%%% fill triangle-packings

\node[vertex, xshift=-17pt] at  (K2botlab) {};
\node[ xshift=-17pt,yshift=-4pt] (v) at  (K2botlab) {$v$};

\begin{scope}[opacity=.5]
\draw[color=cyan,fill=cyan!50] (v) -- ([shift={(\rndcrn+5pt,+8pt)}] K2topsw) -- ([shift={(-\rndcrn-5pt,+8pt)}] K2botne) -- cycle;
\end{scope}

\begin{scope}[opacity=.5]
\draw[color=cyan,fill=cyan!50] (v) -- ([shift={(-5pt,-\crnshift-26pt)}] K1ne) -- ([shift={(-8pt,\crnshift+8pt)}] K1botne) -- 
 ([shift={(8pt,\crnshift+8pt)}] K1botsw) -- cycle;
\end{scope}

\begin{scope}[opacity=.5]
\draw[color=cyan,fill=cyan!50] ([shift={(20pt,8pt)}] K2topsw) -- ([shift={(14pt,-20pt)}] X2sw)
 -- ([shift={(-5pt,-\crnshift-26pt)}] K1ne) -- ([shift={(8pt,\crnshift+-4pt)}] K1botsw) -- cycle;
\end{scope}

%%% draw parts

\draw[thick,rounded corners=\rndcrn] (K1sw) rectangle (K1ne);
\draw[thick,rounded corners=\rndcrn] (K2sw) rectangle (K2ne);

%%% draw division between Ktop and Kbot and separating X

\draw[dashed] (K1topsw)--(K1botne);
\draw[dashed] (K2topsw)--(K2botne);

\draw[color=green!40!black] (X2sw)--(X2se);

%%% draw packing from X

%%% label parts

\draw[thick,color=cyan,fill=white, rounded corners=\rndcrn] ([shift={(5pt,8pt)}] K1topsw) rectangle ([shift={(-5pt,-8pt)}] K1ne);
\node at (K1toplab) {$\Kto$};

\draw[thick,color=cyan,fill=white, rounded corners=\rndcrn] ([shift={(5pt,8pt)}] K2topsw) rectangle ([shift={(-5pt,-8pt)}] X2se);

%\node[thick,rectangle,rounded corners=\rndcrn,draw=cyan, fill=white, inner sep=0pt,minimum width=39pt, minimum height=19pt] at (KX2lab) {};
%\node at (KX2lab) {\tiny $X_2\cap \Ktt$};
%\node[xshift=3pt,yshift=4pt] at (KX2lab) {\tiny $X_1\setminus \Kbt$};
\node[xshift=2.3pt,yshift=3pt] at (KX2lab) {\tiny $X_2\cap \Ktt$};

%\node[color=cyan] at (figlab) {$B$};
\begin{scope}
\def\fromKtx{5.125}
\def\fromKty{-4.35}
\def\fromKtox{3.525}

\coordinate (t1) at ({\fromKtx+.75},\fromKty);
\coordinate (t2) at ({\fromKtx+.50},\fromKty+0.60);
\coordinate (t3) at ({\fromKtx+.25},\fromKty+1.2);

\coordinate (t4) at ({\fromKtox},\fromKty);
\coordinate (t5) at ({\fromKtox},\fromKty+0.60);
\coordinate (t6) at ({\fromKtox},\fromKty+1.20);

\draw[color=cyan, very thick] (t1) -- (t4);
\draw[color=cyan, very thick] (t2) -- (t5);
\draw[color=cyan, very thick] (t3) -- (t6);

%\draw[color=cyan,fill=cyan!20,opacity=0.75] (v) -- (t1) -- (t4)-- cycle;
\draw[color=cyan,opacity=0.75] (v) -- (t1) -- (t4)-- cycle;
\draw[color=cyan,opacity=0.75] (v) -- (t2) -- (t5)-- cycle;
\draw[color=cyan,opacity=0.75] (v) -- (t3) -- (t6)-- cycle;

\node[vertex] at (t1) {};
\node[vertex] at (t2) {};
\node[vertex] at (t3) {};
\node[vertex] at (t4) {};
\node[vertex] at (t5) {};
\node[vertex] at (t6) {};
\end{scope}

\end{tikzpicture}
         \subcaption{set $C'$ of triangles}
         \label{fig:cob-tp2-3}
     \end{subfigure}
        \caption{Triangles in \subref{fig:cob-tp2-1} $A'$, \subref{fig:cob-tp2-2} $B'$, and \subref{fig:cob-tp2-3} $C'$ in the case $x_1> \ell$ and $x_2\ge \ell$.}
        \label{fig:cob-tp2}
\end{figure}
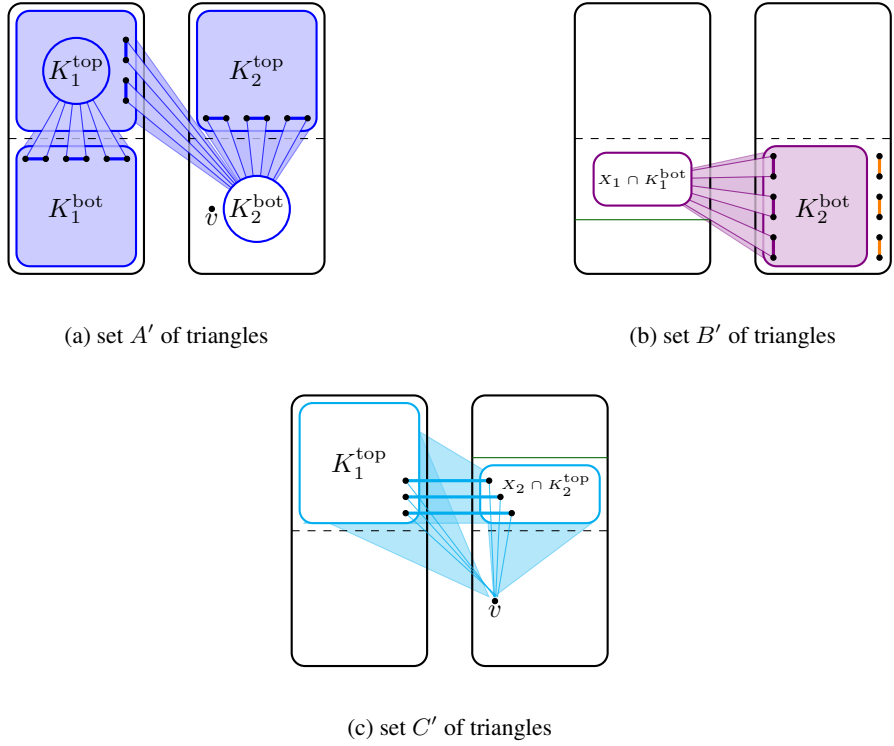
We choose a triangle packing $\TPset$ as follows.
Pack all edges of $\Ktt$ with vertices of $\Kbt$, all edges of $\Kto$ with vertices in $\Kbt$ and all edges of $\Kbo$ with vertices in $\Kto$. 
This gives a set $A'$ of $3\binom{\ell}{2}$ triangles (see~\cref{fig:cob-tp2-1}).
By the second part of \cref{lem:clique_IS} there exists $v \in \Kbt$ such that edges between $v$ and $\Ktt \cup \Kto$ are not used in $A'$.
Additionally, define a set $B'$ of triangles obtained by 
packing edges from $\Kbt$ with vertices of $X_1\cap\Kbo$ (see~\cref{fig:cob-tp2-2}). Then $|B'|=\frac{\ell}{2}(x_1-\ell)$ if $x_1\neq 2\ell$ and $|B'|=\binom{\ell}{2}$ (by~\cref{cor:clique_IS}(\ref{cor:clique_IS:inside})) if $x_1=2\ell$.
Finally, let $C'$ be the set of triangles using $v$ and any maximal matching between $\Kto$ and $X_2 \cap \Ktt$ (see~\cref{fig:cob-tp2-3}).
Since $\Kto$ is complete to $X_2 \cap \Ktt$, we obtain $|C'| = x_2-\ell$. 
It is clear that $\TPset = A' \cup B' \cup C'$ is a triangle packing.

If $x_1<2\ell$, then
\begin{align*}
2\Abs{\TPset} - \Abs{\THset} &\geq 3\ell(\ell-1)+\ell(x_1-\ell)+2(x_2-\ell)-3\ell^2+2\ell - \Brace{x_1 - \ell}\Brace{x_2 - \ell}\\
&= \Brace{x_1 - \ell-1}\Brace{2\ell - x_2} + x_2 - \ell \geq 0.
\end{align*}
The last inequality follows as $x_1\ge\ell+1$.

If $x_1=2\ell$, then we similarly get
\begin{align*}
2\Abs{\TPset} - \Abs{\THset} &\geq 3\ell(\ell-1)+\ell(\ell-1)+2(x_2-\ell)-3\ell^2+2\ell - \ell\Brace{x_2 - \ell}\\
&= \Brace{\ell-2}\Brace{2\ell - x_2} \geq 0.
\end{align*}
We conclude that indeed $2\TP{G}\ge\TH{G}$.
\end{proof}

%    Bibliography.
\bibliographystyle{plainurl}
\bibliography{literature}

\begin{thebibliography}{10}

\bibitem{4colorable}
S.~{Aparna Lakshmanan}, {Cs. Bujt{\'{a}}s}, and {Zs. Tuza}.
\newblock Small edge sets meeting all triangles of a graph.
\newblock {\em Graphs Combin.}, 28(3):381--392, April 2011.
\newblock \href {https://doi.org/10.1007/s00373-011-1048-8}
  {\path{doi:10.1007/s00373-011-1048-8}}.

\bibitem{tightfordensegraphs}
J.~D. {Baron} and J.~{Kahn}.
\newblock {Tuza's conjecture is asymptotically tight for dense graphs}.
\newblock {\em {Comb. Probab. Comput.}}, 25(5):645--667, 2016.
\newblock \href {https://doi.org/10.1017/S0963548316000067}
  {\path{doi:10.1017/S0963548316000067}}.

\bibitem{TuzaEurocomb}
M.~Bonamy, {\L}.~Bożyk, A.~Grzesik, M.~Hatzel, T.~Masařík, J.~Novotná, and
  K.~Okrasa.
\newblock Tuza's conjecture for threshold graphs.
\newblock In J.~Ne{\v{s}}et{\v{r}}il, G.~Perarnau, J.~Ru{\'e}, and O.~Serra,
  editors, {\em Extended Abstracts EuroComb 2021}, pages 765--771, Cham, 2021.
  Springer International Publishing.
\newblock \href {https://doi.org/10.1007/978-3-030-83823-2_122}
  {\path{doi:10.1007/978-3-030-83823-2_122}}.

\bibitem{botler2021tuza}
F.~Botler, C.~G. Fernandes, and J.~Guti{\'e}rrez.
\newblock On {T}uza’s conjecture for triangulations and graphs with small
  treewidth.
\newblock {\em {Discrete Math.}}, 344(4):112281, 2021.
\newblock \href {https://doi.org/10.1016/j.disc.2020.112281}
  {\path{doi:10.1016/j.disc.2020.112281}}.

\bibitem{van2019tight}
W.~{Cames van Batenburg}, T.~{Huynh}, G.~{Joret}, and J.-F. {Raymond}.
\newblock {A tight {E}rd\H{o}s-{P}\'osa function for planar minors}.
\newblock {\em {Adv. Comb.}}, 2019:33, 2019.
\newblock Id/No 2.
\newblock \href {https://doi.org/10.19086/aic.10807}
  {\path{doi:10.19086/aic.10807}}.

\bibitem{chalermsook2020multi}
P.~{Chalermsook}, S.~{Khuller}, P.~{Sukprasert}, and S.~{Uniyal}.
\newblock {Multi-transversals for triangles and the Tuza's conjecture}.
\newblock In {\em Proceedings of SODA 2020}, pages 1955--1974. SIAM, 2020.
\newblock \href {https://doi.org/10.5555/3381089.3381210}
  {\path{doi:10.5555/3381089.3381210}}.

\bibitem{diestelbook}
R.~{Diestel}.
\newblock {\em {Graph theory}}, volume 173.
\newblock Berlin: Springer, 5th edition, 2017.

\bibitem{erdosposa}
P.~{Erd\H{o}s} and L.~{P\'osa}.
\newblock {On independent circuits contained in a graph}.
\newblock {\em {Can. J. Math.}}, 17:347--352, 1965.
\newblock \href {https://doi.org/10.4153/CJM-1965-035-8}
  {\path{doi:10.4153/CJM-1965-035-8}}.

\bibitem{feder2012packing}
T.~Feder and C.~S. Subi.
\newblock Packing edge-disjoint triangles in given graphs.
\newblock {\em Electron. Colloquium Comput. Complex.}, 19:13, 2012.
\newblock URL: \url{http://eccc.hpi-web.de/report/2012/013}.

\bibitem{haxell1999packing}
P.~E. {Haxell}.
\newblock {Packing and covering triangles in graphs}.
\newblock {\em {Discrete Math.}}, 195(1--3):251--254, 1999.
\newblock \href {https://doi.org/10.1016/S0012-365X(98)00183-6}
  {\path{doi:10.1016/S0012-365X(98)00183-6}}.

\bibitem{HaxellRodl}
P.~E. Haxell and V.~R\"{o}dl.
\newblock Integer and fractional packings in dense graphs.
\newblock {\em Combinatorica}, 21(1):13--38, January 2001.
\newblock \href {https://doi.org/10.1007/s004930170003}
  {\path{doi:10.1007/s004930170003}}.

\bibitem{unit}
P.~Heggernes, D.~Meister, and C.~Papadopoulos.
\newblock A new representation of proper interval graphs with an application to
  clique-width.
\newblock {\em Electron. Notes Discret. Math.}, 32:27--34, 2009.
\newblock \href {https://doi.org/10.1016/j.endm.2009.02.005}
  {\path{doi:10.1016/j.endm.2009.02.005}}.

\bibitem{mixed}
J.~Kratochv{\'{\i}}l, T.~Masa{\v{r}}{\'{\i}}k, and J.~Novotn{\'{a}}.
\newblock U-bubble model for mixed unit interval graphs and its applications:
  The {MaxCut} problem revisited.
\newblock {\em Algorithmica}, 83(12):3649--3680, December 2021.
\newblock \href {https://doi.org/10.1007/s00453-021-00837-4}
  {\path{doi:10.1007/s00453-021-00837-4}}.

\bibitem{Krivelevich1995}
M.~Krivelevich.
\newblock On a conjecture of tuza about packing and covering of triangles.
\newblock {\em {Discrete Math.}}, 142(1):281--286, 1995.
\newblock \href {https://doi.org/10.1016/0012-365X(93)00228-W}
  {\path{doi:10.1016/0012-365X(93)00228-W}}.

\bibitem{puleo2015tuza}
G.~J. {Puleo}.
\newblock {Tuza's conjecture for graphs with maximum average degree less than
  7}.
\newblock {\em {Eur. J. Comb.}}, 49:134--152, 2015.
\newblock \href {https://doi.org/10.1016/j.ejc.2015.03.006}
  {\path{doi:10.1016/j.ejc.2015.03.006}}.

\bibitem{JF}
J.-F. Raymond.
\newblock Dynamic {E}rd{\H{o}}s-{P}{\'{o}}sa listing.
\newblock Available at
  \url{https://perso.limos.fr/~jfraymon/Erd\%C5\%91s-P\%C3\%B3sa/}.

\bibitem{tuza1981conj}
{\relax Zs}.~Tuza.
\newblock A conjecture: Finite and infinite sets, {E}ger, {H}ungary 1981, {A}.
  {H}ajnal, {L}. {L}ov{\'a}sz, {V. T. S\'{o}s}.
\newblock In {\em Proc. Colloq. Math. Soc. J. Bolyai}, volume~37, page 888,
  1981.

\end{thebibliography}

\end{document}